\documentclass[11pt]{amsart}
\usepackage{amsmath, amscd, amssymb, latexsym}

\theoremstyle{plain}
\newtheorem{thm}{Theorem}[section]
\newtheorem{lem}[thm]{Lemma}
\newtheorem{cor}[thm]{Corollary}
\newtheorem{pro}[thm]{Proposition}

\theoremstyle{definition}
\newtheorem{df}[thm]{Definition}

\newtheorem{rem}[thm]{Remark}
\newtheorem{ex}[thm]{Example}

\def\om{\omega}
\def\Om{\Omega}

\def\De{\Delta}
\def\ov{\overline}

\def\al{\alpha}
\def\ga{\gamma}
\def\Ga{\Gamma}
\def\si{\sigma}

\def\na{\nabla}

\def\wt{\widetilde}
\def\la{\lambda}
\def\ni{\noindent}
\def\nt{\noindent}

\def\we{\wedge}

\def\lra{\longrightarrow}

\def\be{\beta}
\def\pa{\partial}
\def\th{\theta}
\def\de{\delta}
\def\ti{\tilde}
\def\wt{\widetilde}
\def\mc{\mathcal}

\def\i{{\bf i}}
\def\mr{\mathrm}
\def\R{{\bf R}}
\def\C{{\bf C}}

\def\H{{\bf H}}
\def\Th{\Theta}

\def\b{{\mr{b}}}
\def\f{{\mr{f}}}

\def\bott{\overset{\circ}{\nabla}}
\def\pcon{\overset{\bullet}{\nabla}}

\begin{document}

\title{A higher dimensional foliated Donaldson theory, I}

\author{SHUGUANG WANG}

\address{Department of Mathematics
\\University of Missouri\\Columbia,
MO 65211,USA}

\email{wangs@missouri.edu}

\begin{abstract}
 We introduce a foliated anti-self dual equation for higher dimensional
 smooth manifolds  with codimension-4 Riemannian foliations.
 Several fundamental results including a compactification theorem for the moduli space are formulated and proved,
 towards the defining of a possible Donaldson type invariant for such foliated manifolds.
\end{abstract}

\maketitle

\nt{\em Keywords}: Riemannian foliation, foliated anti-self duality, transverse
anti-self duality, taut foliation, pseudogroup, foliation  cycle.\\

\nt MSC (2010) - Primary: 57R57

\hspace{21mm} Secondary: 53C12,  58H10

\vspace{3mm}

\section{Introduction}\label{int}

 Higher dimensional instanton equations were first proposed by
Physicists in the 1980s, see Corrigan et al \cite{cd}.
Mathematically a program of study was initiated by Donaldson and Thomas \cite{dt} in 1998, where several related
theories were presented, with possible approaches carefully outlined and analyzed.
Since then there  has been a considerable amount of follow-up works by
Thomas \cite{t1},  Donaldson-Segal \cite{ds}, Haydys \cite{fa},
S$\acute{\mr{a}}$ Earp \cite{s}, among others.
Independently Tian \cite{t} laid down the analytical foundation by establishing  the first important
compactification  theorem for the instanton moduli space. A common key feature in  \cite{dt, t, ds} is the natural
occurrences of various ``calibrated'' submanifolds.

It is well known that  gauge theoretic invariants have provided us new and necessary
tools in the understanding  of  three and four dimensional smooth manifolds. The situation is rather
different for higher dimensional manifolds in the sense that, during  the 1970s,
R.C. Kirby and L.C. Siebenmann had already classified smooth structures
on topological manifolds by utilizing Milnor's microbundle theory.
Therefore it would be more interesting and indeed make more sense to place additional geometric structures on
the manifolds and seek possible classifications in the higher dimension case.
Thus \cite{dt}, \cite{dt}, \cite{ds}, \cite{t}, etc. focus on special geometric structures such as Calabi-Yau,
$G_2$ and Spin$(7)$ structures.

In this paper we consider foliation structures and develop an instanton gauge theory for
 smooth manifolds with foliations. By working on this relatively simple situation,
we are able to resolve many of the outstanding issues, as suggested  in the general discussions from \cite{dt, ds,t},
 in order to define the appropriate invariants. Moreover our foliation case provides a good prototype
 to test and confirm a number of  new characteristics about the higher dimensional Donaldson theory.

Specifically, let us take an $n$-dimensional compact smooth manifold $M$, $n\geq 4$, together
 with a transversely oriented Riemannian
foliation $\mc{F}$ of codimension $4$. Given a bundle-like metric $g$, we introduce the foliated anti-self dual
equation for a connection $A$,
$$*(F_A\we \chi)=-F_A,$$
 by invoking the Hodge star $*$ and the $\mc{F}$-characteristic form $\chi=\chi_{\mc F}$
with respect to $g$. Such an equation has in fact appeared in \cite{cd} and \cite{t},
where $\chi$ was taken to  be any differential form of degree $(n-4)$. However in our special
case with foliations, a crucial point made in the paper is that
we apply the equation for basic connections $A$ only. As a consequence of this election,
 the  above equation is equivalent to the transverse
anti-self dual equation $\ov{*}F_A=-F_A$, which is in terms of the Hodge star $\ov{*}$ defined on the
normal bundle $Q=TM/T\mc{F}$ of the foliation. In particular being a foliated anti-self dual connecgtion depends only on the  
component metric $g_Q$ rather
than the whole metric $g$ that is involved in the original equation. This new equation is very much in the sprit of
transverse foliation geometry, which was initiated historically by
 Reinhart \cite{r} and  has been under  extensive research by
A. Haefliger, P. Molino, F.W. Kamber, Ph. Tonduer, to just name  a few. (The other important
way to study foliations is to utilize  A. Connes' non-commutative geometry,
 which is more general and most suitable for the leafwise geometry.)

Regarding the equation above, other than the characteristic form $\chi_{\mc F}$,  one could also take $\chi$ to
be the Kaehler 2-form, the calibration 3-form or the Cayley 4-form respectively,
when $M$ happens to be a complex 3-fold, a $G_2$-manifold
or a Spin(7)-manifold. This is indeed what were explored originally in \cite{dt,ds}.
Since the equation is under-determined, they need to impose additional conditions to secure the ellipticity.
This is comparable to our insistence that the connection $A$ be basic, in case  $\chi$ is
the characteristic form of a foliation.

Following Haefliger \cite{h,h1}, we further make the transition from the equation
$\ov{*}F_A=-F_A$ to a third equation ${*}'F_{A'}=-F_{A'}$, which is now defined on a transversal $Y$
based on the choice of a foliation atlas.
Here the 4-dimensional (non-compact, disconnected) smooth manifold $Y$ comes with a
pseudogroup action and the last equation is the standard four-dimensional anti-self dual equation
that was employed by S.K. Donaldson in his pioneering work during the 1980s.
The issues that $Y$ is non-compact and disconnected
can be overcome since the connection $A'$ on $Y$ is invariant under the pseudogroup action.

In summary, we are guided throughout the paper by the following equivalences:

\hspace{4mm}{Foliated $\chi$-theory on $M$ $\Leftrightarrow$ Transverse theory through  $Q$

\hspace{4mm}$\Leftrightarrow$ Equivariant standard theory on a transversal $Y$.

By using the above equivalences and adapting properly the original Donaldson gauge theory to $Y$,
we  prove a compactification theorem for the foliated instanton moduli space, as a major
improvement of the compactification theorem of Tian \cite{t}  for the more general moduli space.
In particular certain compact leaves of $\mc{F}$ will naturally arise as the branched locus of the limiting
singular connections. With the strategy at hands, we also obtain  complete
results for the linear theory and appropriate perturbations of the foliated moduli space.
The solutions to such issues were only outlined and partially addressed in \cite{dt,ds}.
One key technical assumption required \cite{dt,t,ds} is that the degree $(n-4)$ form $\chi$ be closed.
This is weakened in our foliation case to that the characteristic form $\chi_{\mc F}=\chi$ just needs to
be $\mc{F}$-relatively closed, namely
the foliation $\mc{F}$ is to be taut.

The traditional transverse foliation geometry on $(M, \mc{F})$ developed by Reinhart and others should be viewed as a sort of
geometry on the leaf space $M/\mc{F}$. The point is that $M/\mc{F}$ itself is a highly pathological topological space -- for example
not even Hausdorff in general. Hence it does not allow a direct approach as far as geometry is concerned.
Likewise our foliated Donaldson theory defined on $M$ is really meant to be a sort of Donaldson theory on the 4-dimensional leaf
space $M/\mc{F}$. Clearly then,  when $\dim M=4$, the foliation $\mc{F}$
is trivial and our foliated instanton theory reduces to the original Donaldson theory.

In a sequel, we will consider the orientability issue of the moduli space and further define a Donaldson type invariant
for $(M,\mc{F})$ as a certain multi-linear function on a subgroup of the foliation homology $H^\f_2(M)$.\\

 \ni{Table of contents}

 Section \ref{gen}: Generalities on foliations.

 Section \ref{fse}: Foliated/transverse anti-self dual equations.

 Section \ref{ind}: Transverse ellipticity, index and the virtual dimension.

 Section \ref{pert}: Perturbation of the FASD moduli space.

 Section \ref{comp}: Foliation cycles and compactification of the moduli space.

\section{Generalities on foliations}\label{gen}

Some of the standard references  include  Tondeur \cite{to},
Molino \cite{m}, and Candel-Conlon \cite{cc}. Our purpose  is to
introduce necessary  notations and  results scattering in the literature,
and along the way, to
make  several new observations, which will all be utilized throughout the paper.

To begin, let  $\mc{F}$ denote a $p$-dimensional foliation defined on an $n$-dimensional
compact oriented smooth
manifold $M$. Leaves are then $p$-dimensional immersed submanifolds of $M$.
For the current section, the codimension $q:=n-p$ is arbitrary.

\subsection{Transverse/basic geometry}\label{trab}

In this subsection, no metric is required.
 A foliation atlas
$\mc{U}=\{U_\al,\varphi_\al\}$ yields local foliation coordinates $(x,y)$ with
$x=(x_1,\cdots,x_p)$ and $y=(y_1,\cdots,y_q)$. Thus the last $q$ components of the
coordinate transformation $\ga_{\al\be}$ satisfy
$$\frac{\pa\ga_{\al\be}^{p+j}}{\pa x_i}=0,\; i=1,\cdots,p,\; j=1,\cdots,q.$$
On each chart, the submanifold defined by the equation $y=\mr{const}$ is called a {\em plaque}; plaques can be joined together
to form a leaf.
In terms of the tangent bundle or distribution $L:=T\mc{F}$, the Frobenious theorem
says that $L$ is involutive, i.e. tangential fields  in $\Ga L$ are closed under
Lie bracket. One sets $Q=TM/L$ to be the
normal bundle of $\mc{F}$, yielding the exact sequence
$$0\lra L\lra TM\lra Q\lra 0$$
and its dual $L^*\leftarrow T^*M\leftarrow Q^*\leftarrow0$. In the transverse foliation
theory, $Q$  plays a distinguished role as the
``tangent bundle'' of the leaf space $M/\mc{F}$, which has a highly pathological topology.

By definition, a  {\em foliated principal bundle} $\pi: P\to M$ is a principal bundle
together with a {\em lifted foliation} $\mc{F}^P$ of $\mc{F}$. More precisely at each point $p\in P$ and
its image $z=\pi(p)$,
 $\pi_*$ maps $T_p\mc{F}^P$ isomorphically onto $T_{z}\mc{F}$,
$T_{p}\mc{F}^P\cap T_pP_z=\{0\}$, and $T\mc{F}^P$ is invariant under the structure group $G$-action
on $P$.
 The last two conditions mean that $T\mc{F}^P$ can be viewed as the horizontal distribution of a
{\em partial connection} on $P$.   Since
$T\mc{F}^P$ is involutive, the partial connection is flat. Note also $\dim \mc{F}^P=\dim \mc{F}$.
Thus the lifted foliation $\mc{F}^P$ is such that its leaves are $G$-invariant
and are mapped via $\pi$ diffeomorphically onto leaves of $\mc{F}$.

A {\em foliated vector bundle} $E\to M$ is such that its frame bundle $P_E$ is foliated.
Alternatively and more directly
a foliated vector bundle $E$ is given by a flat {\em partial connection} $\overset{\bullet}{\nabla}$.
Namely $\pcon_Xs$ is defined
for all sections $s\in \Gamma E$ but for tangential fields $X\in \Gamma L$ only, and  it satisfies the
 usual kind of Leibnitz rule
 \begin{equation}\label{lei}
 \pcon_X(fs)=(Xf) s+f\pcon_X s.
 \end{equation}
 Likewise the curvature $\overset{\bullet}{R}(X,X')=\pcon_X\circ\pcon_{X'}-\pcon_{X'}\circ\pcon_{X}-
\pcon_{[X,X']}$ is valid for
tangential fields $X,X'$ only, and is trivial for a flat $\pcon$. A {\em basic section}
$s\in \Ga E$ is  such that
$\pcon_Xs=0,\forall X\in\Ga L$. Because of (\ref{lei}), the concept can be localized: $s$ is
basic iff its restriction
$s|_U$ to any open set  is basic with respect to the partial connection $\pcon|_U$.
One uses $\Ga_\b(E)$ for the set of all basic sections.

Since the lifted foliation on $P_E$ is $GL(r)$-invariant ($r=\mr{Rank}E$), it descends to
a foliation $\mc{F}^E$ on $E$, with leaves mapped diffeomorphically onto the leaves of $\mc{F}$ through
the projection $E\to M$. Then a section $s\in\Ga E$ is basic iff $s$ maps each leaf of $\mc{F}$
diffeomorphically onto a leaf of $\mc{F}^E$.

If $E,E'$ are foliated vector bundles, then the naturally related bundles $E^*, \wedge^rE, \otimes^rE,
 E\oplus E',\mbox{Hom}(E,E')$ etc. all have the
induced foliated bundle structures. One can see this easily
by using the corresponding operations on partial connections. In particular a {\em basic fiber metric}
(Riemannian or Hermitian) $h$ on $E$ is one such that $h\in\Ga_\b(E^*\otimes E^*)$. However
such a basic fiber metric does not always exist in every foliated vector bundle $E$. For example
the existence of a basic metric on $Q$ would mean  that $\mc{F}$ is a {Riemannian
foliation} (see the next subsection).

\begin{ex} The quotient bundle $Q$ is canonically foliated by the Bott
(partial) connection $\overset{\circ}{\nabla}$. This is defined to be
$\bott_X s=\pi{[X,Y]}$, where $X\in \Ga L$ and $Y\in\Ga(TM)$ is a (local) lifting  of
 $s\in \Ga Q$ via the projection  $\pi: TM\to Q$. The Jacobi identity for vector
fields implies the flatness of $\bott$. In contrast, neither $TM$ nor $L$
is a foliated bundle.
 \end{ex}

The Bott partial connection induces one, still denoted by $\bott$, on $\we^r Q^*$. This gives rise
 to  basic sections
 $$\Ga_\b(\we^r Q^*)=\{\al\in \Ga(\we^r Q^*)\mid \bott_X\al=0, \forall X\in\Ga L\}.$$
  These can be viewed as differential
forms on $M$ because of the natural inclusion
 $\Ga(\we^r Q^*)\subset\Ga(\we^r T^*M)=\Om^r(M)$. Indeed one can identify $\Ga_\b(\we^r Q^*)=\Om^r_\b(M)$,
 the latter consisting of {\em basic differential
 forms} $\al\in \Om^r(M)$, that is to say
 \begin{equation}\label{bas}
 \iota_X\al=0, \iota_Xd\al=0, \forall X\in\Ga L
 \end{equation}
where $\iota_X$ is the contraction operator. The first condition means that
$\al$ is reducible to $\Ga(\we^r Q^*)$ while the second means $\al$ is a
basic section on $\we^r Q^*$.
Locally in a foliation chart $(x,y)$, a basic form has the expression
 \begin{equation}\label{bac2}
 \al=\sum_{|J|=r}f_J(y)dy_J
 \end{equation}
where $f_J(y)$ is independent of $x$.
In particular {\em basic functions} are exactly leafwise constant functions.

By the description (\ref{bas}), the exterior differential $d$ preserves basic forms,
giving rise to the basic form complex $\{\Om^*_\b(M), d\}$. The associated De
Rham cohomology $H^r_{\b}(M)$,
$0\leq r\leq q$, is called the {\em basic cohomology} of $\mc{F}$.
One checks easily that
$H^0_{\b}(M)=\R$ and the natural map $H^1_\b(M)\to H^1(M)$ is injective. In general the natural
map $H^r_\b(M)\to H^r(M)$ is neither
injective nor surjective.

Return to the general foliated vector bundle $(E,\pcon)$ together with the flat partial connection.
The basic section set $\Ga_\b(E)$ is a module over basic functions. Under a foliation atlas, transition
functions of $E$ are all basic and vice versa. The bundle $\we^r Q^*\otimes E$ is foliated by
$\bott\otimes\pcon$.  The set of {\em basic forms with values in $E$},
$$\Om^r_\b(E)=\Ga_\b(\we^r Q^*\otimes E)\subset \Om^r(E),$$
is also a module over the space $C^\infty_\b(M)$
of basic functions. Each $\xi\in \Om^r_\b(E)$ can be locally written
as $\xi=\sum\al_i\otimes s_i$
 for some basic forms $\al_i$ and basic sections $s_i$ of $E$.

In the case of foliated principal bundle $P$
with a lifted foliation $\mc{F}^P$, a {\em basic connection} on $P$ is a  regular
connection whose connection form
on $P$ is a basic 1-form with respect to $\mc{F}^P$. (This implies the the horizontal
distribution of the connection contains $T\mc{F}^P$, i.e. the connection is
{\em adapted to} $\mc{F}^P$.)
Then a {\em basic connection} $\nabla$ on $E$  is associated with
a basic connection in the frame bundle of $E$.
Directly on $E$,  a connection $\na$ is basic iff it is
adapted to ${\mc F}$, $\na_X=\pcon_X$,  and its curvature
  satisfies $\iota_X R=0$ for all $X\in\Ga L$.

\begin{rem}\label{neb}
 Although the flat partial connection on  $E$ can always be extended
to an adapted connection, the adapted connection
may not be a basic connection. In fact  $E$  may  not
admit any basic connection at all. The obstruction to the existence
is the Atiyah-Molino class (a secondary characteristic class of the foliated bundle),
see Molino \cite{m} and Kamber-Tondeur  \cite{kt}.
Fundamentally this is due to the lack of a partition of unity $\{\la_\al\}$,
consisting of basic functions, subordinated to any foliation cover $\mc{U}$ of $M$.
(Being basic, $\la_\al$ is constant along each plaque, hence $\mr{Supp}\la_\al\subset U_\al$
can not hold.) It is for the same reason that a foliated vector bundle does not always
admit a basic fiber metric.
\end{rem}

A basic frame $s$ of $E$ is of course
a basic section of the foliated frame bundle $P_E$, i.e. $s\in\Ga_\mr{b}(P_E)$.
Alternatively $s=(s_1, \cdots, s_m)$ where each $s_i$ is a basic section of $E$.
  By restricting the foliation structure
to an open set of $M$, one can also introduce the concept of local basic frames,
which always exist when the open set is small enough.

\begin{lem}\label{eba} For a basic connection  $\na$  in a foliated vector bundle $E$,
the following hold true:

{\em 1)} Locally under a basic frame of $E$, one can write $\na=d+A$ with $A\in\Om^1_\b(\mr{End}E)$.
Consequently the curvature $R$ is a global basic form: $R\in\Om^2_\b(\mr{End}E)$.

{\em 2)} $\na$ preserves basic sections, i.e. $\na:\Ga_\b(E)\to\Om^1_\b(E)$. More generally
the same is true for the extended
differential $d_\na$, namely
$d_\na:\Om^r_\b(E)\to\Om^{r+1}_\b(E)$.

\end{lem}

\nt{\em Proof}. Let $P$ be the foliated frame bundle of $E$. Then $\na$ lifts to a unique
basic connection on $P$; thus the connection 1-form $\om$ on $P$ is basic.

1) A local basic frame $s$ of $E$ comes from a local basic section $t$ of $P$; $t$ being basic
means that it maps local
leaves  of $M$ to local leaves of $P$ diffeomorphically. Then the connection matrix $A$
is the pull-back of $\om$ via $t$, hence $A$ consists of basic form entries since $\om$ is basic.
As for the curvature $R$, locally $R=dA+A\we A$ hence is basic also.

2) Under the local basic frame $s$, each basic section $\xi$ of $E$ is given a vector valued
function $f$ consisting of
basic function components. Hence $\na\xi=(df+Af)s$ is basic. Thus $\na:\Ga_\b(E)\to\Om^1_\b(E)$.
In general $d_\na$ is extended through the Leibnitz rule, and   $d_\na: \Om^r_\b(E)\to\Om^{r+1}_\b(E)$
with $r>0$ follows from the
initial case $r=0$ just shown.
q.e.d.

Note that the lemma is not true if $\na$ is only an adapted connection.
Also we do not have the  operator $\na: \Ga E\to \Ga(Q^*\otimes E)$, despite having
$\na:\Ga_\b(E)\to\Ga_\b(Q^*\otimes E)(=\Om^1_\b(E))$. For instance we have $d:C^\infty_\b(M)\to\Ga_\b(Q^*)$
but not $d: C^\infty(M)\to\Ga(Q^*)$, since $f_y\eta_y dy\not=f_x\xi_y dy+f_y\eta_y dy$
unless $f_x=0$ i.e. $f$ is basic, where $z'=(\xi(x,y),\eta(y))$ is a foliation coordinate change.

It follows from Lemma \ref{eba} that the set $\mc{A}_\b(E)$ of basic connections  forms an affine
space modeled on $\Om^1_\b(\mr{End}E)$. From this and adapting the usual Chern-Weil
transgression argument, one can readily prove the following result.
\begin{pro}\label{bch}
Suppose a foliated complex vector bundle $E$ admits a basic Hermitian metric $h$ and
a compatible unitary basic connection $\na$
with curvature $R$.
Then the class $[\det(1+\frac{\i}{2\pi})R]\in H^*_\b(M)$
is well-defined and is independent of the choices of $\na$ and $h$.
\end{pro}

We will call the resulted individual classes ${c}^\b_k(E)\in H^{2k}_\b(M),
0\leq k\leq \mr{rank}{E}$,
the {\em basic Chern classes} of $E$. Being basic cohomology classes,
${c}^\b_k(E)$ is also subject to $0\leq 2k\leq q$. The basic
Pointryagin classes ${p}^\b_k(E)\in H^{4k}_\b(M)$ of a real foliated bundle
$E$ admitting a basic Riemannian fiber metric and a compatible basic connection
are also defined since the complexified
bundle $E\otimes\C$ carries an induced foliation structure.

\begin{rem}
 In the literature, the characteristic classes of  a foliated bundle belong to the standard
 De Rham cohomology $H^*(M)$, see Kamber-Tondeur \cite{kt, kt1} for example.
  These foliated characteristic classes are  defined by using the
partial flatness of the bundle, hence are  secondary exotic classes.
According to \cite{kt1}, such classes measure
the incompatibility  between the reducibilities of the bundle structure
and the foliation structure into a smaller group.
For instance, in the case of the foliated normal bundle $Q$, of course its
structure group is always reducible to $O(q)$ but its associated Bott
partial connection may not be so. The incompatibility here gives rise to
the Godbillon-Vey class of a foliation.
Foliated characteristic classes are typically  not liftable to the
basic cohomology $H^*_\b(M)$; in fact some of them belong to $H^r(M)$
with $r>q$. Basic characteristic classes and foliated characteristic classes
are best treated as complementary classes, since
the existence of a compatible basic connection as required in  Proposition \ref{bch}
forces most foliated characteristic classes to vanish.
\end{rem}

\subsection{Taut, tight and Riemannian Foliations}\label{RiF}

In this subsection, we will bring  metrics to the geometric
discussions.

Any metric $g$ on $M$ gives
an orthogonal splitting   $TM=L\oplus {L^\perp}$ with $g=g_L\oplus g_{L^\perp}$, and an isomorphism
$L^\perp= Q$.
In  turn $Q$ inherits a metric $g_Q=g_{L^\perp}$. In fact one has the following equivalence
\begin{equation}\label{ge}
g  \Longleftrightarrow (\pi_L, g_L, g_Q)
\end{equation}
where $\pi_L:TM\to L$ is the projection that characterizes the splitting $TM=L\oplus Q$.
Of course $\pi_L$ is subject to $\pi_L^2=\pi_L$. And $\pi_L$ is equivalent to an injection $i:Q\to TM$
such that $\pi\circ i=\mr{Id}$, where $\pi: TM\to Q$ is the quotient projection as before.

Fix a metric $g$ and suppose $\mc{F}$ is tangentially hence transversally oriented.
Choose an oriented orthonormal basis
$$(e_1,\cdots,e_p,\cdots, e_n)\subset L\oplus L^\perp$$
at any point.
The {\em characteristic form} $\chi=\chi_{\mc{F}}\in\Om^p(M)$ of $\mc{F}$ is
defined by
$$\chi(Z_1,\cdots,Z_p)=\mbox{det}(g(e_i,Z_j)_{1\leq i,j\leq p})$$
for any $Z_i\in\Ga(TM)$.
 While
the {\em transverse volume form}
$$\nu\in \Ga(\we^q Q^*)\subset\Om^q(M)$$
is the volume form of $g_Q$ but viewed as a $q$-form on $M$ through
the natural inclusion.

\begin{rem}\label{chi}
It will be useful to note that $\chi$ is the trivial extension of the tangential
volume form of $\mc{F}$ to $TM$ via
the splitting $TM=L\oplus Q$ given by $\pi_L$.
 In terms of the equivalence (\ref{ge}) above,
 $\chi$ depends on $g_L$ and the projection $\pi_L$ alone so is independent of $g_Q$.
On the other hand, $\nu$ depends on $g_Q$ alone.
\end{rem}

Let $\na^M$ denote the Levi-Civita connection on $(TM, g)$. The
mean curvature field (of the leaf submanifolds)  $\tau\in \Ga Q=\Ga L^\perp$ is
defined by $\tau=\sum^n_{i=1}(\na^M_{e_i}e_i)^\perp$.
The mean curvature form is the dual 1-form $\kappa\in\Ga Q^*\subset\Om^1(M)$.

An {\em $\mc{F}$-trivial form} $\varphi$ on $M$ is of degree $\geq p$ and such that
\begin{equation}\label{ftr}
\iota_{X_1}\circ\cdots\circ\iota_{_{X_p}}\varphi=0
\end{equation}
for all tangential fields $X_1,\cdots,X_p\in\Ga L$. A useful fact is the following
Rummler formula \cite{ru}:
\begin{equation}\label{rum}
d\chi=-\kappa\wedge\chi+\varphi
\end{equation}
for some $\mc{F}$-trivial $(p+1)$-form $\varphi$.
This formula is the main reason why $\kappa$ features  prominently in transverse
foliation geometry.

Call $\mc{F}$ is a {\em taut foliation} if there is a metric on $M$ such that
$\kappa=0$, i.e. leaves of $\mc{F}$ are minimal submanifolds.
(Equivalently $\mc{F}$ is {\em harmonic} with respect to this metric, namely the projection
$\pi: TM\to Q$, viewed
as $Q$-valued 1-form in $\Om^1(Q)$, is harmonic.)
 As an application, Rummler \cite{ru,su,h} showed that
$\mc{F}$ is taut iff there is a $p$-form $\xi$ on $M$ such that
$\xi$ restricts non-trivially to $L_z$ at each point $z\in M$ and $\xi$ is
$\mc{F}$-closed (i.e. $d\xi$ is $\mc{F}$-trivial). When $\mc{F}$ is $g$-harmonic,
one can take $\xi=\chi$ to be the $g$-characteristic form.

In a related situation,
$\mc{F}$ is called a {\em tight foliation} if there is a $p$-dimensional calibration on $M$ such that
all leaves are calibrated submanifolds. Then Harvey and Lawson \cite{hl}
showed that
$\mc{F}$ is tight iff there is a closed $p$-form $\xi$ on $M$ such that
$\xi$ restricts non-trivially to $L_z$ at each point $z\in M$. They also showed
that $\mc{F}$ is tight iff $\mc{F}$ is taut and has a complementary distribution that
is integrable.

Taut and tight are for leaf directions. Going transversally,  a metric $g$ on $M$ is called
{\em bundle-like} if the restriction $g_{L^\perp}$  satisfies
 $\mc{L}_X g_{L^\perp}=0$ for all
$X\in \Ga L$. Call $\mc{F}$ a {\em Riemannian foliation} if it admits a bundle-like
metric.

More to the point, $\mc{F}$ being {Riemannian} can be characterized through
a metric $g_Q$ on $Q$ alone without invoking any full metric on $M$. For this
one defines the Lie derivative
\begin{equation}\label{lio}
\mc{L}_X g_Q(s,s')=Xg_Q(s,s')-g_Q(\bott_X s,s')-g_Q(s,\bott_X s')
\end{equation}
for $s,s'\in\Ga Q$ and $X\in\Ga L$, utilizing the Bott partial connection $\bott$ on $Q$.
Then $\mc{F}$ is Riemannian iff $Q$ admits a {\em holonomy invariant metric} $g_Q$,
i.e. $\mc{L}_X g_Q=0$ for all $X\in\Ga L$. The equivalence follows from the fact that
any holonomy invariant metric $g_Q$ on $Q$ extends to a bundle-like metric.
(Indeed, using (\ref{ge}), $g_Q$ together with each metric on $L$ and projection $TM\to L$
determines a unique bundle-like metric.)

For another characterization, a metric $g_Q$ on $Q$ is  holonomy invariant iff
 $g_Q$ is a basic basic metric on $Q$ i.e. $g_Q\in\Ga_\b(Q^*\otimes Q^*)$. The latter is in turn
 equivalent to
 that $g_Q(s,t)$ is a basic function for any basic sections $s,t\in \Ga_\b(Q)$.

The last definition is suitable for a foliated vector bundle $E$:
 one says $E$ is a {\em Riemannian foliated bundle} (or {\em Hermitian foliated bundle})
 if $E$ admits a basic Riemannian (or Hermitian resp.) fiber metric $h$,
 i.e. $h$ is a basic section of $E^*\otimes E^*$.

\begin{ex}\label{rba}
 For a Riemannian foliation $\mc{F}$, the foliated bundle $Q$ admits
a canonical basic connection of some sort.  Indeed given the Levi-Civita
connection $\na^M$ of a bundle-like metric $g$,
the following defines a basic connection $\na'$:
\begin{equation}\label{mer}
\na'_Zs=\left\{\begin{array}{ll}
\bott_Z s\;\; & \mbox{ for } Z\in\Ga L\\
\pi({\na^M_ZY})\;\; & \mbox{ for } Z\in\Ga L^\perp
\end{array}
\right.
\end{equation}
where $Y$ is a local lifting of $s$.  In fact $\na'$ is completely determined
by $g_Q$  and is the only $g_Q$-compatible and torsion free basic connection on $Q$,
see \cite{to}.
\end{ex}

In general, however, a Riemannian (or Hermitian) foliated bundle $E$ does not necessarily
admit any basic connection, whether compatible with the metric or not. This is again because
of the non-existence of basic partition of unity as pointed out in Remark \ref{neb}.
Nonetheless, a Hermitian foliated {holomorphic} bundle over
a complex foliated manifold does admit a basic connection, as one
can check easily the unique connection compatible to both structures is in fact basic.

Being Riemannian is definitely a  restriction for
any foliation $\mc{F}$.
One necessary condition is that all basic cohomology groups
$H^r_\b(M)$ of $\mc{F}$ must be finite dimensional. This is a consequence of
Molino's structure theorem for Riemannian foliations \cite{m}.
In the case of flows (1-dimensional
foliations), any nonsingular Killing field gives a Riemannian flow (which
is also taut).
On the other hand
the famous Anosov flows are all non-Riemannian.

Suppose  now  $\mc{F}$ is Riemannian and fix a bundle-like
metric $g=g_L\oplus g_{L^\perp}$ on $TM=L\oplus L^\perp$. Identify $Q=L^\perp, g_Q=g_{L^\perp}$
as before.
Because $g_Q$ is holonomy invariant, the {\em transverse Hodge star}, $\overline{*}: \wedge^r Q^*\to
\wedge^{q-r} Q^*$,
preserves the basic sections, that is, $\ov{*}:\Ga_\b(\we^r Q^*)\to \Ga_\b(\we^{q-r} Q^*)$. Hence
it passes on to yield
$$\ov{*}:\Om^r_\b(M)\to\Om^{q-r}_\b(M)\;\; \mbox{ for } r=0,1,\cdots,q$$
in view of (\ref{bas}). (Note that $\ov{*}$  will not preserve  basic forms if  $g_Q$  is not
holonomy invariant. Also $\ov{*}$ is not defined in the complement of $\Om^r_\b(M)\subset \Om^r(M)$.)
It relates the usual Hodge star $*:\Om^r(M)\to\Om^{n-r}(M)$ via
\begin{equation}\label{ho}
*\al=\ov{*}\al\we\chi,\; \ov{*}\al=(-1)^{p(q-r)}*(\al\we\chi)
\end{equation}
for $\al\in\Om^r_\b(M)$. (The orientations in $L^\perp, L$ are combined
in such an order to orient $TM$.) In particular  $\chi=*\nu$ but
$\nu=(-1)^{pq}*\chi$.
Also  $\nu$ is a basic form, hence closed. It generates a top degree basic cohomology
class $[\nu]\in H^q_\b(M)$. (Again these use the assumption that $g_Q$ is holonomy invariant.)

Alvarez \cite{a} gave an orthogonal decomposition
$$\Om^r(M)=\Om^r_\b(M)\oplus\Om^r_\b(M)^\perp$$
and showed the basic component $\kappa_\b\in \Om^1_\b(M)$ of $\kappa$ is a closed 1-form.
Hence the twisted exterior differential $d_\kappa=d-\kappa_\b\we$ forms
a complex, resulting the {\em twisted basic cohomology} $\wt{H}^r_\b(M)$.
By using the Hodge decompositions for the  basic Laplacians
$$\De_\b=d_\b\de_\b+\de_\b d_\b,\;
\De_{\b,\kappa}=d_\kappa\de_\kappa+\de_\kappa d_\kappa,$$
 Kamber-Tonduer \cite{kt}
and Park-Richardson \cite{pr} proved the following twisted Poincare duality:

\begin{thm}\label{tpd}
Under a bundle-like metric $g$, the twisted paring
$$\begin{array}{ccl}
\wt{H}^r_\b(M)\times{H}^{q-r}_\b(M)&\lra & {\bf R} \\
([\al],[\be]) & \mapsto & \int_M\al\we\be\we\chi
\end{array}$$
is well-defined and non-degenerate.
\end{thm}

Unlike the classical Poincare duality, the paring here depends on the
 metric $g$ through $\chi$.

Dom{i}nguez \cite{d} has shown that  any Riemannian
foliation $\mc{F}$ carries {\em tense} bundle-like metrics, i.e.
having basic mean curvature form: $\kappa=\kappa_\b$. For such metrics
the twisted cohomology $\wt{H}^r_\b(M)$ is simply from the twisted
differential $d_\kappa=d-\kappa\we$.
To avoid excessive notations with $\kappa_\b$, we will usually use tense bundle-like
metrics in the paper.

\begin{cor}\label{inte}
Suppose $\mc{F}$ is a taut Riemannian foliation so we can choose a bundle-like metric $g$
such that $\kappa=0$. Then the map
\[\begin{array}{ccl}
H^q_\b(M)&\lra& {\bf R} \\

[\alpha] &\mapsto & \int_M\al\we\chi
\end{array}\]
is a well-defined and is an isomorphism. Moreover $[\nu]\in H^q_\b(M)$ is a generator.
\end{cor}

Note that $\kappa=0$ is already needed so that the integral is  well-defined on the cohomology classes.
If the  Riemannian foliation $\mc{F}$ is a non-taut, the integral will depend on the
representatives of each class
and $H^q_\b(M)$ is  trivial, see \cite{t}.

\begin{rem} The main results in Alvarez \cite{a} show that for a Riemannian foliation $\mc{F}$,
 the class $[\kappa_\b]\in H^1_\b(M)$
is independent of the choice of bundle-like metrics. Namely $[\kappa_\b]$ is a topological invariant for
a Riemannian foliation $\mc{F}$. The class $[\kappa_\b]$ vanishes iff $\mc{F}$ is taut.  In particular
any Riemannian foliation on a manifold with $H^1(M)=0$ is taut. On the other hand Anosov
flows are taut but not Riemannian as mentioned above, so Alvarez's theorems do not apply.
\end{rem}


\section{Foliated/transverse anti-self dual equations}\label{fse}

We now take $M$ to be a compact oriented smooth manifold of dimension
$n\geq 4$ and $\mc{F}$ a  transversely oriented Riemannian
foliation of codimension $q=4$.
For convenience we usually work with bundle-like
 metrics that are tense so that the mean curvature
form $\kappa$ is basic hence closed as well.

Consider a foliated $SU(2)$-principal
bundle $P\to M$ and its associated foliated  bundles,
$\mbox{Ad}P=P\times_{\mbox{Ad}}SU(2),
\mbox{ad}P=P\times_{\mbox{ad}}\mathfrak{su}(2), E=P\times_\rho\C^2$, through  the standard (adjoint)
representations
of $SU(2)$. The foliated structures on $\mbox{Ad}P$ and $\mbox{ad}P$ can be best characterized
through the induced flat partial connections.

Let $\mc{A}_\b=\mc{A}_\b(P)$ denote the set of basic connections on $P$.
For $A\in\mc{A}_\b$, consider the usual Yang-Mills functional
$$YM(A)=\int_M |F_A|^2*1=-\int_M \mbox{Tr}(F_A\wedge *F_A),$$
where $F_A$ is the curvature of $A$. (To be more in tune with gauge theory, we use different
 connection and curvature notations  from the last section.)
Since the curvature $F_A$ is a basic form, $F_A\in\Om^2_\b(\mbox{ad}P)$, one
can apply (\ref{ho}) and rewrite the functional in
terms of the transverse Hodge star
\begin{equation}\label{ym}
YM(A)=-\int_M \mbox{Tr}(F_A\wedge \ov{*}F_A)\wedge\chi
\end{equation}where $\chi$ is the characteristic $(n-4)$-form of $\mc{F}$. More generally for
$\al,\be\in\Om^r_\b(\mr{ad}P)$,
\begin{equation}\label{inn}
(\al,\be)=-\int_M \mbox{Tr}(\al \we \ov{*}\be)\wedge\chi
\end{equation}
defines the inner product on $\Om^r_\b(\mr{ad}P)$.

Since $A$ is a basic connection, its extended differential preserves basic forms,
$d_A:\Om^r_\b(\mr{ad}P)\to \Om^{r+1}_\b(\mr{ad}P)$,
by Lemma \ref{eba}.

\begin{pro}\label{tym} The Euler-Lagrange equation of the Yang-Mills functional (\ref{ym})
on $\mc{A}_\b$ is $d^{\ov{*}}_{A,\kappa}F_A=0$, where
\begin{equation}\label{tadj}
d^{\ov{*}}_{A,\kappa}=-\ov{*}(d_A-\kappa\we)\ov{*}
\end{equation}
is the  formal basic adjoint of $d_A:\Om^1_\b(\mr{ad}P)\to \Om^{2}_\b(\mr{ad}P)$.
\end{pro}

\nt{\em Proof}. Let $A\in\mc{A}_\b, \al\in\Om^1_\b(\mr{ad}P)$.
Then $F_{A+t\al}=F_A+td_A\al+t^2\al\we\al$ as usual.
Insert this in the first variation of
(\ref{ym}) and apply Rummler formula (\ref{rum}) for the term $d\chi$ below:
$$\begin{array}{l}
\frac{d}{dt}YM(A+t\al)|_{t=0}\\
\hspace{3mm}=\frac{d}{dt}|_{t=0}\int_M\mr{Tr}(F_{A+t\al}\wedge \ov{*}F_{A+t\al})\wedge\chi\\
\hspace{3mm}=2\int_M\mr{Tr}(d_A\al\we\ov{*}F_A)\we\chi\\
\hspace{3mm}=-2\int_M\mr{Tr}[\al\we (d_A(\ov{*}F_A)\we\chi+ \ov{*}F_A\we d\chi)]\\
\hspace{3mm}=-2\int_M\mr{Tr}[\al\we (d_A(\ov{*}F_A)\we\chi+ \ov{*}F_A\we(-\kappa\wedge\chi+\varphi))]\\
\hspace{3mm}=-2\int_M\mr{Tr}[\al\we (d_A(\ov{*}F_A)\we\chi-\kappa\we \ov{*}F_A\wedge\chi)]\\
\hspace{3mm}=2\int_M\mr{Tr}[\al\we \ov{*} (-\ov{*}d_A(\ov{*}F_A)+ \ov{*}(\kappa\we \ov{*}F_A))\we\chi]\\
\hspace{3mm}=2\int_M\mr{Tr}[\al\we \ov{*}(d^{\ov{*}}_{A,\kappa}F_A)\we\chi]\\
\hspace{3mm}=2(\al,d^{\ov{*}}_{A,\kappa}F_A).
\end{array}$$
Here we drop the term $\al\we\ov{*}F_A\we\varphi$ in the integrand, as it vanishes because $\varphi$ is $\mc{F}$-trivial and
$\al\we\ov{*}F_A$ is basic. (In a local foliation coordinates $(x_1,\cdots x_p, y_1,\cdots, y_q)$,
$\al\we\ov{*}F_A$ contains $d y_j$ only, while $\varphi$ contains no top degree term $dx_1\we\cdots\we dx_p$
in the leafwise $x_i$-portion.)

Set $\frac{d}{dt}YM(A+t\al)|_{t=0}=0$ for every $\al\in\Om^1_\b(\mr{ad}P)$.
This yields the Euler-Lagrange equation $d^{\ov{*}}_{A,\kappa}F_A=0$.
 q.e.d.

Because of the $\kappa$-twisting in $d^{\ov{*}}_{A,\kappa}$,  we will call $d^{\ov{*}}_{A,\kappa}F_A=0$
or explicitly
\begin{equation}\label{tym1}
d_A\ov{*}F_A-\kappa\we\ov{*}F_A=0
\end{equation}
the {\em twisted transverse Yang-Mills equation}.
(Both the usual adjoint $*d_A*$ and the untwisted transverse adjoint
$\ov{*}d_A\ov{*}$ are not suitable here. The first operator does not even  preserve the basic forms.)

If the metric $g$ is not tense, then one needs to replace $\kappa$ with $\kappa_\b$ throughout.
In \cite{gk}, Glazebrook and Kamber  discussed the family index theory of Dirac operators coupled with
basic connections in an auxiliary foliated bundle, generalizing partially one of Atiyah-Singer earlier
works without the presence of foliations.

Note for our codimension-4 case,  the transverse Hodge star
$$\ov{*}: \Om^2_\b(\mr{ad}P)\to\Om^2_\b(\mr{ad}P)$$
satisfies $\ov{*}^2=\mr{Id}$.
So its $(\pm1)$-eigenspaces yields the  decomposition
\begin{equation}\label{stw}
\Om^2_\b(\mr{ad}P)=\Om^+_{\b}(\mr{ad}P)\oplus \Om^-_{\b}(\mr{ad}P)
\end{equation}
which is orthogonal under (\ref{inn}). (One emphasizes that this is a global
decomposition, not pointwise or fiberwise.)
These may be called transverse self-dual and anti-self dual forms respectively.
One can then write $F_A=F^+_A+F^-_A$ for $A\in\mc{A}_\b$.
By definition, $A$ is  {\em transverse anti-self dual} (TASD) if $F_A^+=0$, i.e.
\begin{equation}\label{tasd}
\ov{*}F_A=-F_A.
\end{equation}
In view  of (\ref{ho}), as a crucial fact, the above is equivalent to
\begin{equation}\label{fasd}
*(F_A\we\chi)=-F_A
\end{equation}
utilizing the Hodge star $*$ on $M$ and the characteristic form $\chi$ of ${\mc{F}}$. We will
call (\ref{fasd}) the {\em foliated anti-self dual} (FASD) equation, which is actually defined for all connections $A$.
However when $A$ is basic, we have the equivalence
$$\mr{TASD }\Longleftrightarrow\mr{ FASD}$$
under the bundle-like metric $g$; in particular in this case whether  $A$ is FASD
depends on the holonomy invariant metric $g_Q$ alone, which is perhaps not so
obvious at first glance, since equation (\ref{fasd}) certainly involves the whole metric $g$.

\begin{rem}
Equation (\ref{fasd}) resembles more  the higher dimensional ASD equations
presented in the literature \cite{cd, t}. But one of the main points of our paper is
that by viewing FASD connections as TASD connections, we will be able to apply
directly many results from the classical 4-dimensional gauge theory, particularly when the results
are local in nature.
\end{rem}

Similarly  one can apply (\ref{ho}) and rewrite (\ref{tym1}) so one has the {\em twisted
foliated Yang-Mills equation}
\begin{equation}\label{tfy}
d*(F_A\we\chi)-\kappa\we*(F_A\we\chi)=0.
\end{equation}

In general for an arbitrary Riemannian foliation $\mc{F}$, an FASD connection $A$ may not be
 a minimum for the Yang-Mills functional $YM$.

\begin{pro}\label{topb}
Suppose $\mc{F}$ is a taut Riemannian foliation and choose a bundle-like metric $g$ with $\kappa=0$.
Then the  Yang-Mills functional $YM$ on $\mc{A}_\b$  has a lower bound $8\pi^2\ti{k}$, where
$\ti{k}$ denotes the value $\int_M {c}^\b_2(P)\we\chi$ and ${c}^\b_2(P)$ is the basic Chern class of $P$
as defined in Proposition \ref{bch}.
The lower bound is realized precisely by FASD connections. 
Furthermore $\ti{k}$ depends on the metric $g$, more
precisely it depends on the tangential metric $g_L$ and the projection $\pi_L: TM\to L$
 but is independent of $g_Q$.
\end{pro}

\nt{\em Proof}. As earlier on, consider the foliated vector bundle $E=P\times_\rho {\bf C}^2$ using
the standard representation $\rho:SU(2)\to\mr{Aut}({\bf C}^2)$. Then from Proposition
\ref{bch}, the basic Chern class is
$${c}^\b_2(E)=\frac{1}{8\pi^2}[\mr{Tr}(F_A\we F_A)]\in H^4_\b(M),$$
using any
basic connection $A$ on $P$. (For the induced basic connection on $E$, the curvature remains the same.)
By Corollary \ref{inte}, the integral $\int_M{c}^\b_2(E)\we\chi=\ti{k}$ is well-defined on the
class; in particular the value $8\pi^2\ti{k}$ is independent of $A$.

Recall $F^+_A \perp F_A^-$ and compute formally in
the inner product (\ref{inn}):
$$\begin{array}{ll}
8\pi^2\ti{k}&=\int_M\mr{Tr}(F_A\we F_A)\we\chi\\
&=-(F_A,\ov{*}F_A)\\
&=-(F^+_A+F^-_A, F^+_A-F^+_A)\\
&=-(F^+_A,F^+_A)+(F^-_A,F^-_A)\\
&=YM(A)-2\|F^+_A\|^2.
\end{array}
$$
This shows that $YM(A)\geq8\pi^2\ti{k}$ on $\mc{A}_\b$ and the lower bound is realized by any
FASD connection $A$.

From Remark \ref{chi}, $\chi$ depends on $g_L$ and $\pi_L$ alone; the same is true for
the mean curvature form $\kappa$, see \cite{a} for example. It follows that $\ti{k}$ depends on $g_L$
and $\pi_L$ alone as well.
 q.e.d.

Under the tautness condition, we call the geometry constant
\begin{equation}\label{fchr}
\ti{k}=\ti{k}(P)=\int_M {c}^\b_2(P)\we\chi
\end{equation}
 the {\em foliation charge}
of $P$, in analogy with the classical gauge theory. However here $\ti{k}$  is not quite a topological
 invariant due to its dependence on $g_L$. Nonetheless
we will only perturb the metric $g_Q$ later in Section \ref{pert}, so  the charge $\ti{k}$
will still remain the same.

\begin{rem}
Without the tautness assumption, the integral
$$\int_M \mr{Tr}(F_A\we F_A)\we\chi$$
will depend on the basic connection $A$. (In fact since $H^4_\b(M)=\{0\}$  for any non-taut
Riemannian foliation $\mc{F}$, it seems likely that  there is always  a sequence of basic
connections $A_i$ such that the integral above converges to zero.)
Moreover, an FASD  connection $A$ is not necessarily twisted
 transverse Yang-Mills/twisted
foliated Yang-Mills, as defined in (\ref{tym1}), (\ref{tfy}).
Instead $A$ is  transverse Yang-Mills
 without twisting:
$d_A\ov{*}F_A=0$  (viewing $A$ as TASD in (\ref{tasd})). Applying $d_A$ to
(\ref{fasd}) together with the Bianchi identity, one sees that $A$ also satisfies
\begin{equation}\label{sfym}
d_A*F_A+d\chi\we F_A=0
\end{equation}
which might be called the foliated Yang-Mills equation. Thus an FASD  connection $A$
is also always foliated Yang-Mills.
Regardless of foliations,  equations similar to (\ref{sfym}) have
been proposed in the physics literature \cite{cd}.  However  a
common drawback for all such untwisted higher dimensional Yang-Mills equations is that they do not
arise as the Euler-Lagrange equations of any
 type of energy functionals. Our twisted Yang-Mills equation (\ref{tfy}) is  advantageous
 in this regard, although it is not satisfied by FASD connections.
\end{rem}

\section{Transverse ellipticity, index and the virtual dimension}\label{ind}

Next we take up the linear theory of FASD/TASD, but beginning first with the general theory
on transverse ellipticity.
Suppose $E, E'$ are foliated vector bundles on a manifold $M$ with respect to a foliation $\mc{F}$.

\begin{df}\label{bo}
A {\em basic differential operator} is a linear operator $D:\Ga_\b (E)\to \Ga_\b (E')$
such that locally in
a foliation coordinate $z=(x,y)$ and under basic frames of $E,E'$,
\begin{equation}\label{bdo}
D=\sum_{|\ga|\leq m}C_\ga(y)\frac{\pa^{|\ga|}}{\pa y^{\ga_1}_1\cdots\pa y^{\ga_q}_q}
\end{equation}
where the matrix-valued functions $C_\ga(y)$ all depend on $y$ only. The concept is taken from
\cite{be} and \cite{mo}, which provide some of the background material for this section.
\end{df}

We emphasize that $D$ does not give rise or extend naturally
to a differential operator $\tilde{D}: \Ga E\to \Ga E'$.
 For example a basic connection $\na$ on $E$ induces
a first order basic differential operator
$\na:\Ga_\b (E)\to \Ga_\b (Q^*\otimes E)$ but no natural extension
$\ti{\na}:\Ga(E)\to \Ga(Q^*\otimes E)$ exists; compare with Lemma \ref{eba}
and the remark following it.

Given a basic differential operator $D: \Ga_\b(E)\to\Ga_\b(E')$, one can
define its (transverse) symbol homomorphism $\si:\check{\pi}^* E\to\check{\pi}^* E'$ in the
standard way, where $\check{\pi}:Q^*\to M$ is the projection. More precisely the homomorphism
$\sigma(z,\xi): E_z\to E_z'$
is given by $\sum_{|\ga|= m}C_\ga(y)\xi^\ga$ for any covector $\xi\in Q^*_z\subset T^*_z M$
in the local expression (\ref{bdo}).

\begin{df} $D$ is {\em transverse elliptic}
if its transverse symbol is an isomorphism away from the 0-section, namely
$\sigma(z,\xi)$ is an isomorphism for any $z\in M$ and non-zero $\xi\in Q^*_z$.
\end{df}

Likewise one  defines in the usual way a {\em transverse elliptic complex} of basic differential operators.
One such example is the basic De Rham complex
$$\begin{array}{ccccccc}
&d&&d&&d&\\
\cdots&\rightarrow&\Ga_\b(\wedge^rQ^*)&\rightarrow&\Ga_\b(\wedge^{r+1} Q^*)&\rightarrow&\cdots,
\end{array}$$
whose cohomology yields of course the basic cohomology $H^{\bullet}_\b(M)$.

Suppose $E, E'$ are Riemannian (or Hermitian) foliated vector bundles, so that they admit basic
fiber metrics $h,h'$.
Assume $\mc{F}$ is a Riemannian foliation associated with a bundle-like metric $g$.
Working locally one checks easily that for a basic differential operator $D:\Ga_\b(E)\to\Ga_\b(E')$,
the formal adjoint $D^*:\Ga_\b(E')\to\Ga_\b(E)$ is defined and  basic as well.
(Choose a foliated chart $U\subset M$ with coordinates
$(x,y)$ and basic frames for $E,E'$.
Then $h, h', g, D$ are all represented by (matrix) functions and derivatives in $y$-variables only.
Consider $(D\xi,\eta)=\int_U \eta^T\cdot h\cdot D\xi\cdot\rho dxdy$, where the basic sections $\xi,\eta$,
fiber metric $h$, the derivatives in $D$ and the density $\rho$ of $g$ are all in $y$ only.
Applying the integration by parts gives
$(D\xi,\eta)=(\xi,D^*\eta)$ explicitly, in which $D^*$ surely contains coefficients and differentials
in $y$ only. This verifies $D^*$ is a basic differential operator.) If $D$ is transverse
elliptic then so is $D^*$.

There is no issue to construct various appropriate Sobolev spaces for basic sections $\Ga_\b(E)$.
For example one can restrict the usual Sobolev norms from $\Ga E$ to the subspace
$\Ga_\b(E)$; more directly one can use a basic connection globally or basic functions locally
under basic frames.
We will suppress any explicit Sobolev norms unless it is necessary to specify them otherwise.

By using Molino's structure theorem \cite{m} and a version of equivariant
index, El Kacimi \cite{el} proved the following fundamental result.

\begin{thm}\label{bfi}
Make the afore-stated assumptions. In addition, assume $D$ is transverse elliptic. Then
$D:\Ga_\b(E)\to\Ga_\b(E')$ is Fredholm with a finite index $\mr{Ind}D=\dim\mr{Ker}{D}-\dim\mr{Ker}D^*$.
\end{thm}

\nt{\em  Proof}. We sketch another proof here by working on a transversal manifold $Y$
with a holonomy pseudogroup $\H$ action. The idea is to interpret the
transverse ellipticity of $D$ as the usual ellipticity on $Y$ and apply a
 modified standard pseudodifferential operator theory to
find a parametrix on $Y$. (It will be crucial for us to utilize
the transversal $Y$  in the following sections.)

Choose a (finite) foliation
atlas $\mc{U}=\{U_i,\varphi_i\}$, where $\varphi_i:
U_i\to \R^p\times \R^q$ is a plaque-preserving diffeomorphism, so that $E,E',D,g,h,h'$ etc are all trivialized
appropriately on each chart $U_i$. Assume also each plaque in $U_i$ intersects at most one plaque
in $U_j$ for any $i\not=j$ (shrink the atlas if necessary). Then
the foliation cocycles $\tau_{ij}:\phi_i(U_i\cap U_j)\to\phi_j(U_i\cap U_j)$ can be introduced by matching the
plaques,
where $\phi_i:U_i\to\R^q$ is the composition of $\varphi_i$ with
the projection $\R^p\times \R^q\to\R^q$. Denote by ${\H}$ the holonomy pseudogroup generated by
$\{\tau_{ij}\}$; each element in ${\H}$ is a finite chain of cocycles. For each $i$, fix a
$q$-dimensional submanifold $Y_i\subset M$ so that $\phi_i$ restricts to a diffeomorphism $Y_i\cong\R^q$.
(It will be convenient to simply treat $Y_i=\R^q$ under this diffeomorphism.) The disjoint
 union $Y=\coprod Y_i\subset M$
is a complete transversal of the foliation and is acted upon by ${\H}$ in the sense of a pseudogroup.

By restriction, each $\mc{F}$-basic object in $M$ is in one-to-one correspondence with an
ordinary ${\H}$-invariant object
on $Y$.
For example, a basic differential form on $M$ corresponds to an $\H$-invariant differential
form of the same degree on $Y$.
(To be precise a form $\al$ on $Y$ is $\H$-invariant if $\rho^*(\al|_{R(\rho)})=\al|_{D(\rho)}$ for
any $\rho\in\H$, where $R(\rho), D(\rho)$ are respectively the range and domain of $\rho$.)
 The honolomy invariant metric $g_Q$ on $Q$ associates with an ordinary metric $\ti{g}$ on $Y$ such that
${\H}$ consists of $\ti{g}$-isometries, and vice versa. (Hence $\mc{F}$ is Riemannian iff $Y$ carries an $\H$-invariant
metric.) Likewise the foliated bundle $E$ corresponds
uniquely to a usual bundle $\ti{E}\to Y$ with a lifted ${\H}$-action and the basic sections $\Ga_\b(E)$
correspond to usual ${\H}$-equivariant sections
$\Ga_{\H}(\ti{E})$. The basic differential operator
  $D:\Ga_\b(E)\to\Ga_\b(E')$  of course gives rise to an ${\H}$-equivariant differential operator
$\ti{D}:\Ga(\ti{E})\to\Ga(\ti{E}')$. The transverse ellipticity of $D$ means exactly that $\ti{D}$ is elliptic.
The main idea in our proof is to modify the usual pseudodifferent operator theory to construct a
parametrix for $\ti{D}$. Along the process one must deal with the issue that $Y$ is non-compact,
and this is overcome by utilizing the ${\H}$-invariance. (Implicitly one pretends to work in
the compact but otherwise pathological quotient space $Y/{\H}$.)

Choose a partition of unity $\la=\{\la_i\}$ subordinated to $\mc{U}$. Define a function $\th=(\th_i)$
on $Y$ by
\begin{equation}\label{tla}
{\th}_i(y)=\int_P\la_i\cdot\chi
\end{equation}
where $P=P_y$ is the plaque through $y\in Y_i$ and $\chi$ is the characteristic form of $g$.
Clearly $\th$ has a compact support on $Y$ and is independent of $g_Q$ and
the projection $TM\to L$.

For any $\H$-invariant function ${f}=(f_i)$ on $Y$, introduce its modified (``$\th$-weighted'') integral
\begin{equation}\label{mi}
\int'_Y f d\mu_Y=\sum \int_{Y_i}\th_i(y){f}_i(y)d\mu_{Y_i}
\end{equation}
with the $\ti{g}$-measure $d\mu_Y$. One can use the $\H$-invariance of  ${f}$
to check easily that the integral is actually independent
of the choice of $\la$, though still dependent upon $g_L$.
In the same spirit the modified Fourier transform
$\hat{f}=(\hat{f}_i)$ is defined via
$$\hat{f}_i(\eta)=\int_{Y_i}{\th}_i(y) e^{-\i<y,\eta>}{f}_i(y)d\mu_Y.$$
With the ${\th}$-factor, one further defines the modified pseudodifferential operators
for $\H$-invariant symbol classes from $\ti{E}$ to $\ti{E}'$ as well as the modified
Sobolev spaces. (Recall everything is trivialized on $Y_i=\R^q$.)
 Adjusting the standard proof from Wells \cite{w} for example,
one obtains a parametrix on $\H$-invariant sections for any $\H$-invariant elliptic operator on $Y$.
Consequently and in particular
one has $\dim\mr{Ker}\ti{D}<\infty$ for our  $\ti{D}:\Ga_{\H}(\ti{E})\to\Ga_{\H}(\ti{E}')$.
Since the fiber metrics
$h,h'$ on $E, E'$ are basic, replacing $D$ with its formal adjoint $D^*$, we have
$\dim\mr{Ker}\ti{D}^*<\infty$ for the descendant $\ti{D}^*:\Ga_{\H}(\ti{E}')\to\Ga_{\H}(\ti{E})$.

Pulling back to $M$, we have then that $\dim\mr{Ker}{D}=\dim\mr{Ker}\ti{D}$ and
$\dim\mr{Ker}D^*=\dim\mr{Ker}\ti{D}^*$ are both finite
hence the finiteness of $\mr{Ind}D$ as well.  q.e.d.  \\

Expanding (\ref{mi}), we define the modified integral for any $\H$-invariant top form $\al=(\al_i)$ on $Y$:
$$\int'_Y\al=\sum \int_{Y_i}{\th}_i(y){\al}_i(y).$$
In particular $\int'_Y\ti{\nu}=\mr{vol}_g(M)$ where $\ti{\nu}$ corresponds to the transverse volume
form $\nu$ on $M$.

\begin{rem}\label{patn} The function ${\th}$ defined in (\ref{tla}) plays a vital role here and in later
sections; a few comments are in order.

1)  ${\th}$ is positive on a compact subset of $Y$ that intersects all $\H$ orbits.
If each ${\th}_i$ is extended trivially over other components of $Y$, then the collection
$\{{\th}_i\}$ can be viewed
 as a partition of unity on $Y$ subordinated to the trivial cover $\{Y_i\}$, making the integral definition
 (\ref{mi}) more natural. It is also worth noting that ${\th}$ is defined  whether or not $\mc{F}$
 is Riemannian; indeed it is simply independent of  $g_Q$ as stated above and can be viewed as
 a $\mc{F}$-leafwise datum.
Said differently, the integral (\ref{mi}) on $Y$ incorporates the $\mc{F}$-leaf geometry through ${\th}$.

 2) Haefliger \cite{h} defined a surjective linear map  (``integration along leaves'')
\begin{equation}\label{iaa}
 \int_\mc{F}:\;\;\Om^{p+r}(M)\to\Om^r_c(Y)/\Xi
 \end{equation}
where $\Xi$ is the subspace generated by elements of the
 form $\beta-h^*\beta$ with $\beta\in\Om^r_c(Y),h\in\H$ such that $\mr{supp}\beta\subset\mr{dom}h$.
  As a special case, one computes
 $\int_{\mc{F}}\chi=[{\th}]$ in $\Om^0_c(Y)/\Xi$; thus the class $[{\th}]$
 is  independent of $\la$ just like the integral (\ref{mi}).
It is shown in  \cite{h} that $\mc{F}$ is taut iff $d[\th]=[d\th]=0$ in $\Om^1_c(Y)/\Xi$,
i.e. $[\th]$ might be viewed as a constant function in this case.

For another relevant situation, given a leaf $\ell$ of $\mc{F}$, set
 $$\mr{vol}(\ell)=\sum_i\sum_{y\in\ell\cap Y_i}{\th}_i(y).$$
This defines a leaf volume function, $\mr{vol}: M\to (0,+\infty]$, as studied  in the literature by
D.B.A. Epstein,
R.D. Edwards, K.C. Millet and D. Sullivan \cite{eks}.

 3) By definition, $\th$ is not $\H$-invariant. To compensate for it
one can turn to the ``average'' of $\la_i$ along leaves,
$\ov{\la}_i:M\to [0,1]$, where
$$\ov{\la}_i(z)=\frac{\int_{\ov{\ell}}\la_i d\mu_g}{\int_{\ov{\ell}} d\mu_g}=
\frac{\int_{\ov{\ell}} \la_i d\mu_g}{\mr{vol}(\ov{\ell})}$$
in which $\ov{\ell}$ is the closure of the leaf $\ell$ through $z$.
Because $\mc{F}$ is Riemannian,  $\ov{\ell}\subset M$ is a closed embedded submanifold,
see Chapter 5 of Molino \cite{m}. Each $\ov{\la}_i$ is a basic function as it is already constant
along leaf closures. Thus $\ov{\la}_i$ restricts to an ${\H}$-invariant function $\ti{\la}_i$ on $Y$.
However no $\ti{\la}_i$ is compactly supported, unlike $\th$. Note as well  $\sum\ti{\la}_i=1$ on $Y$.
\end{rem}

Now resume our gauge theory set-up in Section \ref{fse}.
The {\em basic gauge group} $\mc{G}_\b=\Ga_\b(\mr{Ad}P)$ acts on the set
$\mc{A}_\b$, preserving FASD basic connections. Differentiating the action and
linearizing  the equation (\ref{tasd}) at a FASD basic connection
$A$  lead to the fundamental complex
\begin{equation}\begin{array}{ccccc}\label{fec}
&d_A&&d_A^+&\\
\Om^0_\b(\mr{ad}P)&\longrightarrow&\Om^1_\b(\mr{ad}P)&\longrightarrow&\Om^+_{\b}(\mr{ad}P)
\end{array}
\end{equation}
 where $d_A^+$ is the composition of
$d_A$ with the projection onto the basic self-dual $2$-forms.
It is straight forward check that (\ref{fec}) is a transverse elliptic complex.
From Theorem \ref{bfi}, we have the
following  corollary.
\begin{cor}\label{fi} Under a bundle-like metric,
the cohomology groups of (\ref{fec}) have finite ranks so the Euler characteristic ${\chi}^*$ is well defined.
\end{cor}
As usual $-{\chi}^*$ is equal to the index of the ``roll-over'' operator
$d^*_A+d^+_A:\Om^1_\b(\mr{ad}P)\to \Om^0_\b(\mr{ad}P)\oplus \Om^+_{\b}(\mr{ad}P)$.
Write $\ti{d}_P=-\chi^*$ and call it
 the virtual dimension of the FASD moduli space.
Certainly $\ti{d}_P$ is independent of the choice of $A$.

In general it is rather difficult to compute $\mr{Ind}D$
for an arbitrary transverse elliptic basic operator $D$.  Br$\ddot{\mr{u}}$ning el al \cite{bkr} give
an integral formula for $\mr{Ind}D$ which contains some type of eta invariants.
In the special case of gauge theory,
it is tempting to guess that
\begin{equation}\label{vdim}
\ti{d}_P=8 \ti{k}-{3}(1-\ti{b}_1+\ti{b}^+_2)+\zeta
\end{equation}
in terms of  the foliation
charge from (\ref{fchr}) and the basic Betti numbers, where $\zeta$ is some
geometric term depending on the tangential metric $g_L$ only.
Then one must impose that $\mc{F}$ be taut as well.

\begin{rem}
In \cite{f}, H. Fan  considers a somewhat related case where the foliation $\mc{F}$ has one dimension
and the compact manifold
$M$ has non-empty boundary such that  $\pa M$ is transverse to the flows of $\mc{F}$.
A key difference is that his so-called self-duality operator is fully
elliptic on $M$, while our fundamental complex (\ref{fec}) is only {\em transverse} elliptic.
To compute the index, he applies a  classical index formula of P. Gilkey and
evaluates the eta invariant term of $\pa M$ in the formula.  For our more general foliation,
 the relevant eta invariants of \cite{bkr} are defined on some stratifications of singular foliations
 which are not directly tied to a non-compact transversal
$Y$ of the foliation.

In the special case of Sasakian manifolds (among Riemannian taut
flows), Biswas and Schumacher \cite{bs} introduced and discussed the moduli spaces
of Sasakian stable bundles.
\end{rem}

%
%

\section{Perturbation of the FASD moduli space}\label{pert}

Still use the set up in Section \ref{fse}. Call a basic connection $A\in\mc{A}_\b$ {\em reducible}
if $\mr{ad}P=L'\oplus L''$ for some foliated line
bundles $L', L''$ and  $A=A'\oplus A''$ correspondingly for some basic connections
$A', A''$.

\begin{pro}\label{irr}
Suppose $A\in\mc{A}_\b$ is not flat. Then the following are equivalent:

{\em 1)} $d_A:\Om^0_\b(\mr{ad}P)\to
\Om^1_\b(\mr{ad}P)$ has a non-trivial kernel;

{\em 2)} $A$ is reducible;

{\em 3)}  $\mc{G}_{\b,A}/{\pm1}$ is non-trivial, where $\mc{G}_{\b,A}$ is the
stabilizer of the $\mc{G}_{\b}$ action at $A$.

\end{pro}

\nt{\em Proof}. This models verbatim the original proof in Freed-Uhlenbeck \cite{fu}.

1) $\Rightarrow$ 2). Take any non-trivial $u\in\mr{Ker}d_A$. Then as in \cite{fu}, $u$ has a constant
$\i$-eigenvalue $\la$ and eigenvector $e$, $ue=\la e$. Because $u$ is basic,
$e$ is basic as well. Set $L'=<e>, L''=(L')^{-1}$. Since $d_A e=0$ and $A$ is basic,
$A=A'\oplus A''$ splits correspondingly into basic connections.

2) $\Rightarrow$ 3). Under $\mr{ad}P=L'\oplus L''$, the acting group
$\mc{G}_{\b,A}/{\pm1}$ contains a
subgroup $\mr{diag}(e^{\i\th},e^{-\i\th})$.

3) $\Rightarrow$ 1). Take any non-trivial element $u$ in the Lie algebra of
$\mc{G}_{\b,A}/{\pm1}$. Then  $d_A u=0$ and of course $u$ is basic.  q.e.d.  \\

Thus $\mc{G}_{\b}/{\pm1}$ acts freely on the set of irreducible basic connections $\mc{A}_\b^*$.\\

As a base metric, let  us fix a bundle-like but not necessarily tense metric $g=g_L+g_{L^\perp}=g_L+g_Q$, with
the induced
identifications $TM=L\oplus L^{\perp}=L\oplus Q$. Unlike the usual 4-dimensional
gauge theory, we can not take all perturbations in $\mc{N}=\Ga(\mr{GL}(TM))$ since they
do not preserve bundle-like metrics. Instead we use only those in $\mc{N}_\b=\Ga_\b(\mr{GL}(Q))\subset\mc{N}$.
Clearly for each $\varphi\in\mc{N}_\b$, $\varphi^*g$ is also bundle-like and at the same time
has the same tangential component $g_L$. Plainly put,  the basic automorphisms $\mc{N}_\b=\Ga_\b(\mr{GL}(Q))$
parameterize the set of all holonomy invariant metrics on $Q$ while
fixing the tangential metric. We will use $\mc{N}_\b$ as the perturbation space for our
moduli space of FASD connections.

\begin{thm}\label{tray} Let $P_+:\Om^2_\b(\mr{ad}P)\to \Om^+_\b(\mr{ad}P)$ be the projection
onto the $g_Q$-self dual basic forms.
Then the map
$$\begin{array}{clcr}
\mc{P}=\mc{P}_\b: &\mc{A}^*_\b\times\mc{N}_\b&\lra&\Om^+_\b(\mr{ad}P)\\
&(A,\varphi)&\mapsto&P_+((\varphi^{-1})^*F_A)
\end{array}$$
has 0 as a regular value. For ease on the notations, we have suppressed the  $C^k$ and Sobolev $L^2_{l}$ norms
that must be imposed on the spaces above.
\end{thm}

\nt{\em Proof}. As alluded in the proof of Theorem \ref{bfi}, the key
 will be to use the  transversal manifold $Y$ introduced in that proof.
Since $Y$ is 4-dimensional,  we have the usual pointwise self-duality
on $Y$. Note that in the original proof of Theorem 3.4 in
Freed-Uhlenbeck \cite{fu}, most of the arguments
are pointwise in nature.
One can modify and adapt their proof for $Y$, keeping in mind that our
$Y$ is non-compact, disconnected as well as with a pseudogroup action.

Let $(A,\varphi)$ be a zero of $\mc{P}$ and we need to show that
the differential
$$\de \mc{P}=\de \mc{P}_{A,\varphi}: T_A(\mc{A}^*_\b)\oplus T_\varphi(\mc{N}_\b)\lra
\Om^+_\b(\mr{ad}P)$$ is surjective at $(A,\varphi)$, namely $\mr{Coker}(\de \mc{P})=0$.
The bundle-like metric $g$ defines the $L^2$-inner
product $(\;,\;)$ on $\Om^+_\b(\mr{ad}P)$. In turn
\begin{equation}\label{ck}
\Phi\in\mr{Coker}(\de \mc{P})\Leftrightarrow
 (\de \mc{P}(\xi), \Phi)=0
\end{equation}
 for all $\xi\in\mr{Dom}(\de \mc{P})=\Om^1_\b(\mr{ad}P)\oplus
 \Ga_\b(\mathfrak{g}\mathfrak{l}(Q))$.

Of course the foliated bundle $P$ corresponds uniquely to
an  $SU(2)$-bundle $P'$ on $Y$ with a lifted $\H$-action, where $\H$ is the
holonomy pseudogroup defined in the proof of Theorem \ref{bfi}.
In the same spirit, $g_Q$ corresponds uniquely to an $\H$-invariant metric $g'_Q$ on $Y$
and basic sections in $\Om^2_\b(\mr{ad}P)$ to $\H$-invariant sections
in $\Om^2_{\H}(\mr{ad}P')$. From $\mc{P}(A,\varphi)=0$,
$A$ corresponds to an $\H$-invariant {{ASD connection}} $A'$ on $P'$ with respect to
$\varphi'^*(g'_Q)$, where $\varphi'\in\Ga(\mr{GL}(TY))$  corresponds to $\varphi$.
Moreover $\mc{P}$ corresponds to the obviously defined $\mc{P}'$.

To show $\de \mc{P}$ is surjective, it is equivalent to showing that
$$\de\mc{P}'= \mc{P}'_{A',\varphi'}: T_{A'}(\mc{A}'^*_\H)\oplus T_{\varphi'}(\mc{N}'_\H)
\lra \Om^+_\H(\mr{ad}P')$$
is surjective, where $\mc{N}'_\H=\Ga_\H(\mr{GL}(TY))$. To have the appropriate regularity
property,
one must use suitable Sobolev spaces for the above spaces, which are constructed in the proof
of Theorem \ref{bfi}.
Pick any $\Phi'\in\mr{Coker}(\de\mc{P}')$ and we need to show $\Phi'=0$. From
 the transverse elliptic complex (\ref{fec}) and its corresponding {\em elliptic}  complex on $Y$,
and proceeding as in \cite{fu}, one shows that $\Phi'$ is a function of sufficient
smoothness rather than just a distribution.

 Let $\Phi\in\mr{Coker}(\de \mc{P})$ correspond
to $\Phi'$; $\Phi$ is a function also. From (\ref{ck}) and using (\ref{tadj}),
one has the pointwise equation
\begin{equation}\label{uts}
\ov{*}d_A\ov{*}\wt{\Phi}-\ov{*}(\kappa_\b\we \wt{\Phi})=0,
\end{equation}
where $\wt{\Phi}=\varphi^*(\Phi)$, $\ov{*}$ is with respect to $\varphi^*g_Q$ and
$\kappa$ is the mean curvature form of $\varphi^*g$.
The basic component $\kappa_\b$ is used since
$\varphi^*g$ is not necessarily tense. Varying the second component of $\de\mc{P}$ one also
has the pointwise equation
\begin{equation}\label{ap}
(r^*F_A,\wt{\Phi})_{\varphi^*g}=0
\end{equation}
for all $r\in\Ga_\b(\mathfrak{g}\mathfrak{l}(Q))$.
Because everything is basic, equations (\ref{uts}) and  (\ref{ap}) correspond respectively to the pointwise
equations on $Y$
\begin{equation}\label{sm}
{*}'d_{A'}{*}'{\wt{\Phi}}'-{*}'(\kappa'_\b\we \wt{\Phi}')=0,\; (r'^*F_{A'},\wt{\Phi}')_{\varphi'^*g'_Q}=0
\end{equation}
where $*'$ is the $\varphi'^*(g'_Q)$-Hodge star and $r' \in\Ga_\H(\mathfrak{g}\mathfrak{l}(TY))$.
Applying Lemma 3.7 of \cite{fu} to the second equation of (\ref{sm}), $\mr{Im}(F_{A'})$ and
$\mr{Im}(\wt{\Phi}')$ are pointwise perpendicular; hence one of $F_{A'}, \wt{\Phi}'$ has
rank 1 whenever they are both non-zero.

We now move  entirely to the 4-dimensional manifold $Y$ with the metric $\varphi'^*(g'_Q)$.
Since $A'$ is ASD, together with the usual Bianchi identity, $F_{A'}$ satisfies the elliptic linear
equation
$$\Delta_{A'}F_{A'}=0,$$
 where $\Delta_{A'}=(d_{A'}^*d_{A'}+d_{A'}d_{A'}^*)$ and $d_{A'}^*=-*'d_{A'}*'$.
 By definition, $\wt{\Phi}'$ is self-dual. Together with the first equation in (\ref{sm}),
 this implies that $\wt{\Phi}'$ also satisfies an elliptic linear equation:
$$\Delta_{A'}\wt{\Phi}'-d_{A'}^*(\kappa'_\b\we \wt{\Phi}')+d_{A'}{*}'(\kappa'_\b\we \wt{\Phi}')=0.$$
As a consequence of these elliptic equations,  a slight generalization of the unique continuation theorem
says that if $F_{A'}$ or $\wt{\Phi}'$
vanishes in some open set of $Y$ then it vanishes identically. (Here is a quick proof the
generalization. From the construction
in the proof of \ref{bfi},  $Y=\coprod_{i=1}^m Y_i$,
each $Y_i$ being  an open and connected. Suppose $F_{A'}$ vanishes on some open set $U\subset Y$
and $U\cap Y_1\not=\emptyset$. By the original Aronszajin's unique continuation theorem,
$F_{A'}$ vanishes entirely on $Y_1$. By the $\H$-invariance,  $F_{A'}$ vanishes on
the open set $\mr{Im}(\tau_{12})\subset Y_2$, where  $\tau_{12}$ is the foliation cocycle
used in the proof of \ref{bfi}. By Aronszajin's theorem again, $F_{A'}$ vanishes on
$Y_2$. Continue the same argument, then $F_{A'}$ vanishes on the whole $Y$.)

Following \cite{fu} roughly, the rest of the proof proceeds like this.
Suppose on the contrary that $\Phi'\not=0$ so $\wt{\Phi}'=\varphi'^*(\Phi')\not=0$ at some
point. From the unique continuation theorem above, $\wt{\Phi}'$ does not vanish
in some open and dense subset of $Y$.
We will derive a contradiction in both of the following cases:

\nt Case 1. $\mr{Rk}(F_{A'})=2$ at some point in $Y$ hence in an open neighborhood.
From the previous and an earlier paragraphs, $\wt{\Phi}'$ has rank 1 in a possibly smaller
neighborhood $Z$. Set $d_{A',\kappa}=d_{A'}-\kappa'_\b\we$. Equation (\ref{sm}) and the self-duality
of $\wt{\Phi}'$ implies $d_{A',\kappa}\wt{\Phi}'=0$.
Write $\wt{\Phi}'=\al\otimes u$ with $\al$ a self-dual
2-form and $u$ a unit section of $\mr{ad}P$ in some open set of $Z$. Because $\kappa'_\b$ is
of length 1, $d_{A',\kappa}\al\otimes u=0$ implies $d_{A',\kappa}u=0$ similar to \cite{fu}.
On the other hand, since $\kappa_\b'$ is 1-form, $d_{A',\kappa}^2u=d_{A'}^2u=[F_{A'},u]$,
which is non-zero by Lemma 3.12 of \cite{fu}. This contradicts the just derived $d_{A',\kappa}u=0$  .

\nt Case 2. $\mr{Rk}(F_{A'})\leq 1$ on $Y$. Since $A'$ is not flat, the above unique continuation
theorem implies that $F_{A'}$ has {rank} 1 in some open and dense subset $V=\{F_{A'}=0\}^c$ of $Y$.
Locally in $V$, write $F_{A'}=\si\otimes u$ as before,  with $\| u\|=1$.
From $d_{A'}F_{A'}=0$ (Bianchi identity, no $\kappa$-twist anymore), one gets $d_{A'}u=0$.
Using the same argument as \cite{fu}, one shows that $V\cap Y_i$ is connected for each $i$.
Since $d_{A'}u=0$, $u$ extends over $V^c=\{F_{A'}=0\}$ to $\hat{u}$ on each $Y_i$.
Therefore $d_{A'}\hat{u}=0$ and $\hat{u}\not=0$, contradicting that $A$ namely $A'$
is irreducible.

Hence we must have $\Phi'=0$ which means $\Phi=0$ and consequently $\mr{Coker}\de\mc{P}=0$.
That is to say, $0$ is a regular value of $\mc{P}$.  q.e.d.

The standard argument brings us the following perturbation theorem.

\begin{thm}\label{pmo}
The parameterized moduli space $\ov{\mc{M}}^*_b:=\mc{P}_\b^{-1}(0)/\mc{G}_\b\subset(\mc{A}^*_\b/\mc{G}_\b)\times\mc{N}_\b$
is a smooth Banach manifold. The projection map $\ov{\pi}_\b:\ov{\mc{M}}^*_\b\to \mc{N}_\b$ is Fredholm
of index $\ti{d}_P$. Consequently there is a Baire subset $\{\varphi\}$ of $\mc{N}_\b$,
whose FASD moduli spaces $\mc{M}^*_\b(\varphi^*g)$ are all smooth of dimension $\ti{d}_P$.
\end{thm}

\nt{\em Proof}. Here the $\mc{G}_\b$ action on $\mc{A}^*_\b$ has a slice given
by $\mr{Ker} d^{\ov{*}}_{A,\kappa_\b}$, in terms of the twisted formal adjoint as in Proposition \ref{tym}.
Together with \ref{tray}, the tangent space of $\ov{\mc{M}}^*_b$ can then be represented
by $\mr{Ker}(\de\mc{P})\cap \mr{Ker} d^{\ov{*}}_{A,\kappa_\b}$.
That $\ov{\pi}_\b$ is Fredholm follows from
$$\mr{Ker}(\de\ov{\pi}_\b)\cong \H^1_\b, \mr{Coker}(\de\ov{\pi}_\b)\cong\H^+_\b$$
where $\H^1_\b,\H^+_\b$ are the cohomology groups of the complex (\ref{fec}). This
also identifies the index with that of (\ref{fec}), since $\H^0_\b$ is trivial because
of irreducibility and Proposition \ref{irr}. Finally the standard Sard-Smale theorem yields
the Baire set as the regular value set of $\ov{\pi}_\b$.  q.e.d.

\begin{rem}\label{tper}
The perturbation theorem can be adapted easily to  taut foliations. In this case, choose
the base metric $g$ so that $\mc{F}$ is a harmonic Riemannian foliation.
 Recall the mean curvature form $\kappa$ depends only on the tangential
metric component $g_L$ and  the projection $\pi:TM\to L$ (a consequence of
Rummler's formula (\ref{rum}), see also \cite{a}),
 hence it will remain to be trivial for all the perturbed
 metrics $\varphi^*g$, $\varphi\in\mc{N}_\b$. In particular for those generic $\varphi$
 the FASD moduli spaces $\mc{M}^*_\b(\varphi^*g)$ are all smooth and come with the same
 foliation charge $\ti{k}$.
 \end{rem}

%
%

\section{Foliation cycles and compactification of the moduli space}\label{comp}

Denote by $\mc{T}_r=\mc{T}_r(M)$ the set of all $r$-dimensional currents on $M$, namely
the set of topological linear functionals on $\Om^r_c(M)=\Om^r(M)$.
Let $\mc{F}^{r}$ denote the set of $\mc{F}$-trivial $r$-forms, $r\geq p$,
as defined in (\ref{ftr}). For convenience, also set $\mc{F}^{r}=\Om^r(M)$ when $r<p$.
The annihilator set is $\mr{Ann}(\mc{F}^r)=\{T\in\mc{T}_r: T \mbox{ vanishes on } \mc{F}^r\}$.
Inspired by Haefliger \cite{h}, we propose the following definition.
\begin{df}\label{fcr}
 $T\in\mc{T}_r$ is called a {\em foliation current}
if $T\in\mr{Ann}(\mc{F}^r)$
and $\pa T\in\mr{Ann}(\mc{F}^{r-1})$.
\end{df}
Use $\mc{T}_r^\f$ for the set of all $r$-dimensional foliation currents.
Note that $\mc{T}_0^\f=\cdots=\mc{T}_{p-1}^\f=\{0\}$ and any $p$-dimensional
foliation current  is automatically closed, since
$$\mr{Ann}(\mc{F}^{0})=\cdots=\mr{Ann}(\mc{F}^{p-1})=\{0\}.$$
The resulted {\em foliation homology}
$$H^\f_r(M)=\mr{Ker}\{\pa:\mc{T}_r^\f\to \mc{T}_{r-1}^\f\}
/\mr{Im}\{\pa:\mc{T}_{r+1}^\f\to \mc{T}_{r}^\f\}$$
satisfies $H^\f_0(M)=\cdots= H^\f_{p-1}(M)=\{0\}$. These trivial groups present a valid
reason to  truncate and re-set ${H}^\f_r(M):=H^\f_{p+r}(M)$,
where $r=0,\cdots, q$.

Recall each differential form $\tau\in \Om^r(M)$ defines a diffuse
current $\tau^\we\in\mc{T}_{n-r}$ through $\tau^\we(\om)=\int_M\tau\we\om$.
If $\tau$ is a basic form, it is not hard to confirm that $\tau^\we$ will
be a foliation current. Furthermore Sergiescu \cite{se} proves the following
De Rham duality theorem for the basic cohomology.

\begin{thm}\label{ser}
For any Riemannian foliation $\mc{F}$, there is a canonical isomorphism
$$H^r_\b(M)\xrightarrow{\cong} {H}^\f_{q-r}(M)$$
induced by the diffuse map $[\tau]\mapsto[\tau^\we]$.
\end{thm}

Now specialize to $p$-dimensional currents. Even though a foliation $p$-current
$T$ is always closed,
$T$ is not quite a {foliation cycle} according to our next definition.

\begin{df}\label{fcc} A {\em foliation cycle} is a foliation $p$-current
$T \in \mc{T}_p^\f$
that is represented by integration (i.e. a current of integration).
\end{df}
This means that $T$ can be extended to the set of continuous $p$-forms on $M$ as a
continuous functional.
 By adding this condition, one can characterize foliation cycles in several
other ways, which appear quite distinctive at first glance.

\begin{pro}\label{fce}
The following are equivalent statements for any $p$-dimensional current $T$:

{\em 1)}  $T$ is a foliation cycle.

{\em 2)} $T$ is closed and lies in the topological closure of the set of
finite linear combinations of Dirac foliation currents.

{\em 3)} $T$ is a closed current of integration whose tangent simple $p$-vector
field is linearly dependent with  that of $\mc{F}$, namely
$\overrightarrow{T}=\lambda\overrightarrow{\mc{F}}\in\Ga(\we^p(T\mc{F}))$ a.e.
for some scalar function $\la$.

{\em 4)} $T$ corresponds to a holonomy invariant signed (Borel) measure $\mu$
on a transversal $Y$ via ``integration along leaves''.
\end{pro}

In part 2), a Dirac foliation current is one given by a $p$-vector in
the leaf tangent space $T_x\mc{F}$ at a point $x\in M$.
The various equivalences in Proposition \ref{fce} are proved respectively by
 Ruelle-Sullivan, Harvey-Lawson and Haefliger:
 2) $\Leftrightarrow$ 4) in \cite{rs} (see also Theorem 10.2.12 of \cite{cc}),
 2) $\Leftrightarrow$ 3) in \cite{hl},
 and 1) $\Leftrightarrow$ 4) in \cite{h}.
One can also impose that the
foliation current $T$ in 1) be {\em positive}, which corresponds to $\la>0$ in 2),
 positive coefficients in the linear combinations in 3), and a positive measure $\mu$
 in 4).

The correspondence 1) $\Leftrightarrow$ 4) is particularly relevant for us,
for it serves as a prelude to the proof of our convergence theorem \ref{conv} below.
Explicitly, given a holonomy invariant (i.e. $\H$-invariant) signed measure $\mu$ on $Y$,
the corresponding foliation cycle $T_\mu:\Om^p(M)\to\R$ is given by
\begin{equation}\label{ial}
T_\mu(\om)=\sum_i\int_{Y_i}\left(\int_{P_y}\la_i\om\right)d\mu
\end{equation}
where $\{\la_i\}$ is a partition of unity subordinated to the
foliation cover $\{U_i\}$ and $P_y$ is the plaque through $y$.
More generally, given a 0-dimensional $\H$-invariant current $T'$ on $Y$
(not necessarily a measure), the analogous formula
\begin{equation}\label{iay}
T(\om)=\sum_i T'\left(\int_{P_y}\la_i\om\right)
\end{equation}
establishes a one-to-one correspondence between
 $\H$-invariant $0$-currents on $Y$ and foliation $p$-currents on $M$.
In fact for $\H$-invariant  currents $T'$ of arbitrary dimension $r$ on $Y$,
the same formula (\ref{iay}) still makes sense and gives again an isomorphism with
$(p+r)$-dimensional foliation currents $T$ on $M$.
(Both (\ref{ial}) and (\ref{iay}) are independent of the choice of $\{\la_i\}$.)
The last isomorphism further appears as the bottom map in the following commutative diagram:
\begin{equation}\label{dc}
\begin{array}{ccc}
\Om^{q-r}_\b(M)&\xrightarrow{\cong}&\Om^{q-r}_\H(Y)\\
\downarrow& &\downarrow\\
\mc{T}_{p+r}^\f(M)&\xleftarrow{\cong} &\mc{T}_r^\H(Y).
\end{array}
\end{equation}
Here $\Om^{q-r}_\H(Y), \mc{T}_r^\H(Y)$ contain $\H$-invariant forms and currents respectively on $Y$.
The top horizontal map is the restriction map which is also an isomorphism. The vertical maps are the injective
diffuse homomorphisms (the left one induces the isomorphism in Theorem \ref{ser}).

Moreover it is shown  in \cite{h} that $\mc{T}_r^\H(Y)= (\Om^r_c(Y)/\Xi)^*$
 and the bottom isomorphism  is the dual of the linear map $\int_\mc{F}$
stated in (\ref{iaa}). (Namely $\mc{T}_{p+r}^\f(M)\subset(\Om^{p+r}(M))^*$ is the image set of
the dual map $\int^*_\mc{F}$, which is of course injective.) Our definition \ref{fcr} for foliation currents
 is  motivated by this isomorphism and Theorem \ref{ser}.
Originally foliation currents were defined differently in \cite{rs,hl} to be those currents,
but not necessarily closed,
 from parts 2), 3) of Proposition \ref{fce}. Their definitions are valid
for dimension $p$ currents only, while \ref{fcr} in our case is applicable for all dimensions.

Recall a rectifiable current $T$ is the limit, under the mass norm, of Lipschitz images  of finite
polyhedral chains with integer coefficients, see \cite{ff}. Equivalently $T$ is given by
its the support $S$ (a rectifiable subset of $M$) and integer multiplicity $\Th$ through the formula
$T(\om)=\int_S \om(\overrightarrow{S})\Th d\mc{H}^r$, where $\overrightarrow{S}$ is
the tangent simple $p$-vector field of $S$ and $\mc{H}^r$ is the Hausdorff measure.
Here is our main convergence theorem on the FASD moduli space $\mc{M}_b(E)$.

\begin{thm}\label{conv} Suppose $\mc{F}$ is a codimension-4 taut transversely oriented Riemannian foliation
on a compact oriented smooth manifold $M$ of dimension $\geq 4$.
Let $\{A_\al\}$ be a sequence of FASD basic connections defined on a foliated $SU(2)$-bundle $E\to M$.
Then there exits a rectifiable foliation cycle $T=(S,\Th)$ and an FASD connection $\ti{A}$
defined in another foliated $SU(2)$-bundle $\ti{E}$ such that the following hold:

{\em 1)} There is a subsequence $\{A_{\beta}\}$ converging to $\ti{A}$ on $M\backslash S$ after
applying basic bundle isomorphisms $\rho_{\beta}: E|_{M\backslash S}\to \ti{E}|_{M\backslash S}$.

{\em 2)} As foliation cycles, the diffuse currents $\frac{1}{8\pi^2}(\mr{Tr}( F_{A_{\beta}}^2))^\we$
 converge weakly  to $\frac{1}{8\pi^2}(\mr{Tr}(F_{\ti{A}}^2))^\we +T$.

 {\em 3)} The foliation charges of $E,\ti{E}$ are constrained by
 $$\ti{k}(E)=\ti{k}(\ti{E})+{\bf M}(T)$$
 where ${\bf M}$ is the mass norm for currents.

{\em 4)} The support $S$ of $T$ consists of a finite number of compact leaves,
$S=\{L_1,\cdots, L_w\}$. With the multiplicity $\Th=(\th_1,\cdots,\th_w)$,
the mass norm ${\bf M}(T)=\th_1\mr{vol}(L_1)+\cdots+\th_w\mr{vol}(L_w)$,
thus depending on the tangential metric $g_L$ only.
\end{thm}

\nt{\em Proof}. We continue working with the transversal $Y=\coprod^m_{i=1} Y_i$ and holonomy
pseudogroup $\H$. Unless otherwise said,  notations will be
carried over from the proofs of Theorems \ref{bfi} and \ref{tray}.
We further choose the foliation atlas $\mc{U}=\{U_i,\varphi_i\}_{i=1}^m$
to be {\em regular} in the sense of Definition 1.2.11 of Candel-Conlon \cite{cc}
(see their Lemma 1.2.17 for the existence).
Essentially this means that the closure $\ov{U}_i$ is compact and is
contained  in another  foliation chart $W_i$ as foliated sets.
Hence the plaque ${P}_y$ of $U_i$ through any $y\in Y_i$ has the closure $\ov{P_y}$
 contained in a plaque of $W_i$ and the extended plaque $P_y$ of $U_i$
is defined even if $y\in \ov{Y_i}$.

The foliated $SU(2)$-bundle $E$ associates a unique $SU(2)$-bundle $E'\to Y$ with a lifted
$\H$-action. Since $\mc{F}$ is Riemannian, the transverse volume form $\nu$ is basic
so corresponds to an $\H$-invariant top form $\nu'$ on $Y$. Of course $\nu'$ is
simply the volume form of $g'_Q$.

By Proposition \ref{topb} and since $A_\al$ is FASD, i.e. TASD, one has for each foliation
chart $U_i$,
$$\begin{array}{ll}
8\pi^2\ti{k}&=\int_M\mr{Tr}(F_{A_\al}\we F_{A_\al})\we\chi=\int_M|F_{A_\al}|^2\nu\we\chi\\
&\geq\int_{\ov{U}_i}|F_{A_\al}|^2\nu\we\chi=\int_{\ov{Y_i}}(\int_{\ov{P_y}}\chi)|F_{A'_\al}|^2\nu',
\end{array}$$
where $A'_\al$, corresponding to $A_\al$, is $\H$-invariant ASD  on $E'\to Y$.
In the above we have applied the Fubini theorem for the iterated integral and the fact that
$|F_{A_\al}|^2$ is a basic function. Set
$$K_i=\min_{y\in \ov{Y}_i}\int_{\ov{P}_y}\chi.$$
Since the restricted $\chi$ is the $g_L$-volume form of $\ov{P_y}$ and each  $\ov{P_y}$ is
contained in a plaque of $W_i$,
one sees that $K_i$ is strictly positive (and finite). It follows that
$$\int_{Y_i}|F_{A'_\al}|^2\nu'\leq 8\pi^2\ti{k}/K_i$$
or more uniformly
\begin{equation}\label{xc}
\int_{Y}|F_{A'_\al}|^2\nu'\leq 8\pi^2\ti{k}/K
\end{equation}
if one sets $K=\min\{ K_1,\cdots,K_m\}$. From its definition, $K$ depends on the tangential metric $g_L$.

By (\ref{xc}), the measures $\frac{1}{8\pi^2}|F_{A'_\al}|^2\nu'$ are bounded, hence  weakly converge  to a
measure: $\frac{1}{8\pi^2}|F_{A'_\al}|^2\nu'\to\mu$ on $Y$ (taking a subsequence if necessary). Since
$\int_Yd\mu\leq \ti{k}/K$, for any $\epsilon>0$, there are at most $\ti{k}/(\epsilon^2K)$
``blow-up'' points $\{y_1,\cdots, y_l\}$ in $Y$ which do not lie in
any $g'_Q$-geodesic ball of $\mu$-measure
less than $\epsilon^2$. Hence around each point $y\in \Om:=Y\backslash\{y_1,\cdots, y_l\}$,
there is a small geodesic
ball $D_y\subset Y$ such that $|F_{A'_\al}|^2\nu'\leq\epsilon^2$ for large $\al$. Choose $\epsilon$\
small enough so that a version of Uhlenbeck's theorem, Proposition (4.4.9) of Donaldson-Kronheimer
\cite{dk}, applies to give us a (sub)sequence of gauge transformations $u'_\al$ on $E'|_\Om$
with $u'_\al(A'_\al)$ converging to $A'$ on $\Om$. Of course
$\int_{Y}|F_{A'}|^2\nu'\leq 8\pi^2\ti{k}/K$. From the 4-dimensional Uhlenbeck removable singularity
theorem (Theorem 4.4.12 of \cite{dk}), one extends $A'$ smoothly to another $g'_Q$-ASD connection
$\ti{A}'$ on $\ti{E}'\to Y$, where the restriction  $\ti{E}'|_\Om$ is isomorphic to $E'$. By
checking the proof
of Uhlenbeck's theorem, one sees that $\ti{E}'$ carries a lifted $\H$-action  and $\ti{A}'$
is $\H$-invariant. Thus they correspond to a foliated bundle $\ti{E}$ and an FASD-basic connection
$\ti{A}$ respectively on $M$.

From $\frac{1}{8\pi^2}|F_{A'_\al}|^2\nu'\to \mu$  on $Y$ and
$|F_{A'_\al}|^2=|F_{u'_{\al}(A'_\al)}|^2\to |F_{A'}|^2$ on $\Om$,
one can write
\begin{equation}\label{mbn}
\mu=\frac{1}{8\pi^2}|F_{A'}|^2\nu'+\sum_{r=1}^l n_r\de_{y_r}
\end{equation}
where $n_r=\mu(f_r)-\frac{1}{8\pi^2}\int_Y f_r|F_{A'}|^2\nu'$ and $f_r$ is any compact-supported function
that equals $1$ near $y_r$ and $0$ near each $y_j\not=y_r$. Choose a small ball $B$ centered at $y_r$
such that $f_r|_B\equiv 1$. Then
$$\begin{array}{ll}
8\pi^2 n_r&=\lim_{\al\to\infty}\int_B f_r|F_{A_\al}|^2\nu'-\int_B f_r|F_{A'}|^2\nu'\\
&=\lim_{\al\to\infty}\int_B\mr{Tr}(F_{A_\al}\we F_{A_\al})-\int_B\mr{Tr}(F_{A'}\we F_A)\\
&=\lim_{\al\to\infty}\int_B(\mr{Tr}(F_{u_\al(A_\al)}\we F_{u_\al(A_\al)})-\mr{Tr}(F_{A'}\we F_{A'}))\\
&\equiv\lim_{\al\to\infty}\int_{\partial B}(\mr{CS}(u_\al(A_\al))-\mr{CS}(A'))
\;\;\;\mr{mod }\;\;8\pi^2\mathbb{Z}\\
&=0
\end{array}$$
where $\mr{CS}$ stands for the Chern-Simons form. Thus $n_r$ is an integer, as in the case of the
classical 4-dimensional gauge theory.

Let $\hat{S}$ be the set of leaves each of which passes through at least one  point
in $W:=\{y_1, \cdots, y_l\}$. Define $S=\cup\{L:L\in\hat{S}\}$ to be the underlying
set of points in $\hat{S}$, i.e. ${S}$ is the saturation of $W$.
Because  $|F_{A_\al}|^2$ are basic functions on $M$,  each element in
$S\cap Y$ is a blow-up point for the sequence
$\{|F_{A'_\al}|^2\}$.
Thus $S\cap Y=W$ hence a finite set. Consequently each leaf  $L$ in $\hat{S}$
must be compact; otherwise, in view that  $Y$ is a complete transversal,
$L$ being non-compact would imply that $L$ intersects $Y$ infinitely many times.
(More precisely note that
each foliation chart $U_i$ has a compact closure and $L\cap U_i$
contains $|L\cap Y_i|$  plaques,  which is a finite number less than  $|S\cap Y|$.
It follows from Theorem 5 on page 51 of \cite{cn} that $L$ must be closed.)

Since the measure $\mu$ is $\H$-invariant, the set $W$ is also $\H$-invariant.
Write $\hat{S}=\{L_1,\cdots, L_w\}$ explicitly. Then $y_r$ is $\H$-equivariant to $y_{s}$ iff
$y_r$ and $y_{s}$ belong to a common leaf $L_j\in \hat{S}$. Said differently
each leaf leaf $L_j\in \hat{S}$ represents an $\H$-equivalence class in $W$
and $w$ is then the number of elements in the quotient set $W/\H$.

From (\ref{mbn}), the set $\{n_1,\cdots, n_l\}$ is also $\H$-invariant, since the
 other two terms in
the equation are so. This means that $n_r=n_{s}$ if $y_r$ and $y_{s}$ are
 $\H$-equivariant. Thus the common integer $n_r=n_{s}$ can be associated with
 the leaf $L_j\in \hat{S}$ and  will be denoted by ${\th}_j$, regarded as
 the multiplicity of $L_j$.
Pooling together we have a $p$-dimensional current $T$ with support $S$ and multiplicity
$\Th=(\th_1,\cdots,\th_w)$,
$$T=L_1+\cdots+ L_1+\cdots+ L_w+\cdots+ L_w,$$
where each $L_j$ is repeated $\th_j$ times.

Certainly $T$ is a rectifiable current. By condition 3) in Proposition \ref{fce},
$T$ is a foliation cycle. Alternatively  we check that $T$ is associated
with the $\H$-invariant Dirac measure $\ov{\de}=\sum^l_{r=1} n_r\delta_{y_r}$  through the
formula (\ref{ial}): for any $\om\in\Om^p(M)$,
$$\begin{array}{ll}
T(\om)&=\sum^w_{j=1}\th_j\int_{L_j}\om\\
&=\sum^w_{j=1}\th_j\int_{L_j}\sum_{i=1}^m\la_i\om\\
&=\sum_i\sum_j\th_j\int_{L_j}\la_i\om\\
&=\sum_i\sum_j\th_j\int_{L_j\cap Y_i}\la_i\om \;\; (\mbox{since supp}\la_i\subset U_i)\\
&=\sum_i\sum_j\th_j\sum_{y_r\in L_j\cap Y_i}\int_{P_{y_r}}\la_i\om\\
&=\sum_i\sum_j n_r\sum_{y_r\in L_j\cap Y_i}\int_{P_{y_r}}\la_i\om\;\;
(\mbox{from the definition of }\th_j) \\
&=\sum_r n_r\int_{P_{y_r}}\la_i\om=\sum_r n_r\de_{y_r}(\int_{P_y}\la_i\om).
\end{array}$$
That is to say, $T(\om)=\sum_i\int_{Y_i}( \int_{P_{y}}\la_i\om)d\ov{\de}$.
(Note the function $\int_{P_y}\la_i\om$ does  have compact support in $y$.)
By 4) of Proposition \ref{fce}, $T$ is a foliation cycle again. (To be logically
correct,  (\ref{ial}) also requires that the  foliation cover $\{U_i\}$ be regular.)

Now turn to calculate
$$\begin{array}{ll}
&\int_{M}\frac{1}{8\pi^2}\mr{Tr}( F_{A_\al}^2)\we\om
=\frac{1}{8\pi^2}\int_{M}|F_{A_\al}|^2\nu\we\om\\
&\;\;=\frac{1}{8\pi^2}\int_{M}|F_{A_\al}|^2\nu\we(\sum_i\la_i)\om)\\
&\;\;=\frac{1}{8\pi^2}\sum_i\int_{Y_i}(\int_{P_y}\la_i\om)|F_{A'_\al}|^2\nu'\;\;(\mbox{Fubini's theorem})
\end{array}$$
 From the weak convergence $\frac{1}{8\pi^2}|F_{A'_\al}|^2\nu'\to\mu$ and
 (\ref{mbn}), as $\al\to\infty$, the last term converges to
 $$\frac{1}{8\pi^2}\sum_i\int_{Y_i}\left(\int_{P_y}\la_i\om\right)
|F_{{A}'}|^2\nu'+\sum_r n_r\de_{y_r}\left(\int_{P_y}\la_i\om\right)$$
which equals $\int_{M}\frac{1}{8\pi^2}\mr{Tr}( F_{\ti{A}}^2)\we\om+T(\om)$ as we have just exhibited.
Thus we have shown the diffuse current sequence
$\frac{1}{8\pi^2}(\mr{Tr}(F_{A_\al}^2))^\we$ converges weakly to $
\frac{1}{8\pi^2}(\mr{Tr}(F_{\ti{A}}^2))^\we +T$.
 Since $\mr{Tr}F_{A_\al}^2$ is a basic form, $\frac{1}{8\pi^2}(\mr{Tr}(F_{A_\al}^2))^\we$
is a foliation current hence a foliation cycle as well (it is a current of integration
for sure). The same is true for $\frac{1}{8\pi^2}(\mr{Tr}(F_{\ti{A}}^2))^\we$.
This completes the proof of part 2) of the theorem.

From the weak convergence in 2), we have, as $\al\to\infty$,
$$\int_{M}\frac{1}{8\pi^2}\mr{Tr}( F_{A_\al}^2)\we\chi\to
\int_{M}\frac{1}{8\pi^2}\mr{Tr}( F_{\ti{A}}^2)\we\chi+T(\chi),$$
which leads to the charge equation in 3).

Part 4) and the rest of the theorem are pretty clear.  q.e.d.

\begin{rem}
 1) In essence, Theorem 4.2.3 of Tian \cite{t} implies that the current $T$
 appearing in our Theorem \ref{conv} is
 a rectifiable foliation  cycle. Then Corollary 1.16 of \cite{hl}
shows also that the support $S$ of $T$ consists of a finite number of
compact leaves. Our proof here is much shorter and Theorem
\ref{conv} contains a lot stronger results than \cite{t}. For example the limit
connection $\ti{A}$
extends smoothly on the entire manifold $M$, meaning the singular set
$S([\ti{A}])$ is empty. In the general case of \cite{t}, it is only conjectured that
$S([\ti{A}])$ has Hausdorff codimenison $\geq 6$ (see p. 263 therein).

2) The tautness assumption in \ref{conv} is weaker than $\chi$ being
closed, which is as assumed in \cite{t} (cf. the discussion on tightness in Subsection \ref{RiF}).
Even without the tautness assumption, our arguments above show
that parts 1), 2), 4) of Theorem  \ref{conv} continue to hold,
as long as the sequence $\int_{M}\frac{1}{8\pi^2}\mr{Tr}( F_{A_\al}^2)\we\chi$
is bounded for all $\al$. The essential fact needed here is that a basic connection $A$ on $E$
is FASD iff the corresponding $\H$-invariant connection $A'$ on $E'$ is ASD, regardless whether
$\mc{F}$ is taut. (In contrast, $A$ is {\em twisted} foliated Yang-Mills (\ref{tfy})
iff $A'$ is Yang-Mills without the tautness.)
\end{rem}

The existence of a compact leaf is a fundamental problem in foliation theory.
On one extreme, a theorem of Molino \cite{m} says that for a Riemannian foliation
$\mc{F}$ with all leaves compact, the leaf space $M/\mc{F}$ is an orbifold.
The converse is evidently true as well.

To state the compactification theorem, introduce the set
 $\Upsilon_E=\{(\hat{E}, T)\}$ where $\hat{E}$  is a foliated $SU(2)$-vector bundle
and  ${T}=\{L_1,\cdots,L_n\}$ consists of compact leaves, repeated according to multiplicity.
The pair $(\hat{E},S)$ is  subject to the condition
\begin{equation}\label{cstt}
\ti{k}(\hat{E})+\mr{vol}({L})=\ti{k}({E})
\end{equation}
where $\mr{vol}({T})=\mr{vol}({L}_1)+\cdots+\mr{vol}({L}_n)$.
 Define the ideal FASD moduli
space
$$\mc{I}\mc{M}_\b(E)=\{([A],T)\mid [A]\in\mc{M}_\b(\hat{E}), (\hat{E}, T)\in \Upsilon_E\},$$
whose topology is characterized through the sequence convergence: a sequence
$([A_\al],T_\al)$ in $\mc{I}\mc{M}_\b(E)$ is convergent to $([A],T)$ if after suitable basic
gauge transformations $A_\al \to A$
on $M\backslash (\cup_\al\mr{supp}T_\al\cup \mr{Supp}T)$ and
$T_\al\to T$ weakly as currents.

\begin{thm}\label{comt} Suppose $\mc{F}$ is a taut foliation and $g$ is a bundle-like
metric with mean curvature form $\kappa=0$. Then the ideal moduli space $\mc{I}\mc{M}_\b(E)$ is compact.
Moreover $\mc{I}\mc{M}_\b(E)$ depends on the leafwise metric component of $g$ alone.
\end{thm}

\nt{\em Proof}. The compactness is a consequence of Theorem \ref{conv}. The sole dependence on $g_L$
is due to that Equation (\ref{cstt})  involves $g_L$ only in the foliation charge and leaf volumes. q.e.d.
\\

In a future paper, we plan to address the orientability issue of the moduli space ${\mc M}_\b={\mc M}_\b(E)$. We will use a natural
foliation structure on $\mc{B}$ to define
a slant product $H^4_\b(\mc{B}\times M)\times H^\f_2(M)\to H^2_\b(\mc{B})$. Then we define a Donaldson invariant
as a multi-linear functional on $H^\f_2(M)$ by making a suitable paring with
the fundamental class $[{\mc M}_\b]$.
In order to get an integer-valued invariant we restrict to a subgroup of $H^\f_2(M)$,
whose classes are generated by rectifiable currents. In the case of all leaves being compact, we would also like to
compare our foliated theory with the existing orbifold Donaldson theory, which is defined on the orbifold leaf space $M/\mc{F}$.

A similar strategy is also being pursued  for
Seiberg-Witten invariants by A. Renner and J. Lee.

\end{document}